\documentclass[twoside,11pt]{article}

\usepackage{blindtext}

\usepackage{jmlr2e}
\usepackage{amsmath}
\usepackage{enumerate}
\usepackage{url} 
\usepackage{amssymb}
\usepackage{bm}
\usepackage{subfigure}
\usepackage{mathrsfs}
\usepackage{bbm}
\usepackage{color}
\usepackage{algorithmicx,algorithm}
\usepackage{epstopdf}
\usepackage{cases}
\DeclareMathOperator*{\argmin}{arg\,min}

\def\tr{{\rm tr}}
\def\prox{{\rm Prox}}
\def\TgFDR{{\rm TgFDR}}

\def\A{{\mathcal A}}
\def\B{{\mathcal B}}

\def\D{{\mathcal D}}
\def\E{{\mathcal E}}
\def\M{{\mathcal M}}
\def\N{{\mathcal N}}
\def\O{{\mathcal O}}
\def\Y{{\mathcal Y}}
\newcommand{\blue}[1]{\begin{color}{black}#1\end{color}}

\usepackage{lastpage}

\ShortHeadings{Penalized Low-Rank Tensor Regression}{Chen and Luo}
\firstpageno{1}

\begin{document}

\title{Group SLOPE Penalized Low-Rank Tensor Regression}

\author{\name Yang Chen\thanks{Corresponding author.} \email yangchen@bjtu.edu.cn \\
       \addr School of Mathematics and Statistics\\
       Beijing Jiaotong University\\
       Beijing, P. R. China
       \AND
       \name Ziyan Luo \email zyluo@bjtu.edu.cn \\
       \addr School of Mathematics and Statistics\\
       Beijing Jiaotong University\\
       Beijing, P. R. China}

\editor{My editor}

\maketitle

\begin{abstract}
This article aims to seek a selection and estimation procedure for a class of tensor regression problems with multivariate covariates and matrix responses, which can provide theoretical guarantees for model selection in finite samples. Considering the frontal slice sparsity and low-rankness inherited in the coefficient tensor, we formulate the regression procedure as a group SLOPE penalized low-rank tensor optimization problem based on an orthogonal decomposition, namely TgSLOPE. This procedure provably controls the newly introduced tensor group false discovery rate (TgFDR), provided that the predictor matrix is column-orthogonal. Moreover, we establish the asymptotically minimax convergence with respect to the \blue{TgSLOPE estimate risk}. For efficient problem resolution, we equivalently transform the TgSLOPE problem into a difference-of-convex (DC) program with the level-coercive objective function. This allows us to solve the reformulation problem of TgSLOPE by an efficient proximal DC algorithm (DCA) with global convergence. Numerical studies conducted on synthetic data and a real human brain connection data illustrate the efficacy of the proposed TgSLOPE estimation procedure.
\end{abstract}


\section{Introduction}
Tensor regression modeling, in which the regression coefficients take the form of a multi-way array or tensor, is an important and prevalent technique for coefficient estimation and/or feature selection in the high-dimensional statistical learning theory, with wide applications in many modern data science problems such as neuroimaging analysis, see, e.g., \cite{2013CP,2017STORE,2019Convex,2020Tucker,2022AOS}. In this article, we focus on the tensor regression with multivariate covariates and matrix responses. Given $n$ independent and identically distributed (i.i.d.) observations $\{(\bm{x}_i,\bm{Y}_i)\}_{i=1}^n$ with $\bm{x}_i\in\mathbb{R}^{p}$ the vector of predictors and $\bm{Y}_i\in\mathbb{R}^{p_1\times p_2}$ the matrix of responses, the regression model can be expressed as follows
\begin{eqnarray}\label{tm-1}
\bm{Y}_i=\B^*\times_{3}\bm{x}_i+\bm{E}_i, i\in[n]:=\{1,2,\dots, n\},
\end{eqnarray}
where $\B^*\in\mathbb{R}^{p_1\times p_2\times p}$ is the unknown coefficient tensor with some inheritance structures like sparsity and/or low-rankness, $\bm{E}_i, i=1,\dots,n$, are i.i.d. noise matrices whose entries are i.i.d. drawn from the Gaussian distribution $N(0,\sigma^2)$. Our goal is to seek a feature selection and tensor estimation procedure for model (\ref{tm-1}). A straightforward idea to estimate $\B^*$ is via optimization problems
\begin{equation*}
\hat{\B}=\argmin_{\B\in\Omega} f(\B; \D),
\end{equation*}
where $\Omega$ is any constraint set of sparse and/or low-rank tensors, $f(\B; \D)$ can be taken as the least squares loss or any more general loss function with $\D=\{(\bm{x}_i,\bm{Y}_i)\}_{i=1}^n$ the random sample set.

\subsection{Related Work}
One simple approach to estimate the sparse and/or low-rank coefficient tensor is to use matricization techniques such that the model (\ref{tm-1}) is degenerated into a linear matrix regression model, in which the estimated tensor $\B^*$ and the response matrix $\bm{Y}_i$ are unfolded into a $p \times p_1p_2$ matrix and a $p_1p_2$ dimensional vector, respectively. Therefore, \blue{a huge body of relevant work on sparse and low-rank matrix methods} can be applied to deal with the tensor regression problems, see, e.g., variable selection by sparsity penalized methods \citep{OWJ:2011,2019Convex,JMVA2021}; low-rank matrix estimation by reduced-rank regression methods \citep{Lowrank2018,Nuclear2019,2021SRRR} \blue{and nuclear norm penalized methods \citep{LRMM-20, LRMM-21}}; \blue{joint penalized methods for selection and low-rank estimation \citep{LRMM-22-1}; regularized covariance estimation for matrix-valued data \citep{LRMM-22-2}}. However, the use of matricization techniques will not only break the sparse and low-rank structures of tensors, making resulting estimators difficult to interpret, but also lead to a dramatic increase in dimensionality, which is prone to over-fitting phenomenon.

Keeping the tensor format for regression models, the existing work can be divided into two categories from the perspective of different characterizations for the sparsity and low-rankness of coefficient tensors. One is based on the tensor decomposition, including CP decomposition \blue{(CPD)} \citep{2013CP,2017STORE,2020Cubic} and Tucker decomposition \citep{2018Tucker,2020ISLET,2022AOS,2022JCGS}. The other is the method without tensor decomposition, such as imposing sparsity on elements, fiber vectors, and slice matrices \citep{2019Convex}, assuming low-rankness of slice matrices \citep{2019Non,2019Convex,2020L2RM}, and considering Tucker low-rankness of coefficient tensors \citep{2019Non}.

Focusing on the matrix response tensor regression model (\ref{tm-1}), limited work has been done on statistical property analysis and algorithm design. \cite{2017STORE} have considered an element-wise sparse tensor regression model based on the \blue{CPD}, and established a non-asymptotic estimation error bound for the estimator obtained from the proposed alternating updating algorithm, which compounds the truncation-based sparse tensor decomposition procedure. Considering the frontal slice sparsity and low-rankness, \cite{2020L2RM} have proposed a two-step screening and estimation procedure and shown that it enjoys estimation consistency and rank consistency. \cite{2022JCGS} have developed a generalized Tucker decomposition model with features on multiple modes, and investigated the statistical convergence for the proposed supervised tensor decomposition algorithm with side information. Nevertheless, these methods may not work well for finite samples, since their feature selection results are usually achieved in infinite samples, or unfortunately some of these are not capable of feature selection.

To seek a mechanism that enables to make inference about the validity of the selected model in finite samples, \cite{2015slope} have introduced a new convex penalized method for classical linear regressions inspired by the Benjamini-Hochberg procedure \citep{1995BH}, namely Sorted L-One Penalized Estimation (SLOPE). They have shown that SLOPE controls the false discovery rate (FDR) under certain conditions. Following \cite{2015slope}, a recent line of research \citep{2016slope,2018slopeAOS,2019gslope,2019slopeJMLR,2021SRRR} has studied the SLOPE based methods, including statistical property analysis, algorithm design, etc. One of these work \citep{2019gslope} has extended SLOPE to group SLOPE (gSLOPE) to deal with the situation when one aims to select whole groups of regressors instead of single one. This motivates us to investigate the use of gSLOPE penalty for feature selection and estimation in the tensor regression framework with finite samples.

\subsection{Our Contributions}
We consider the matrix response tensor regression model (\ref{tm-1}) and embed the frontal slice sparse and low-rank structures for the estimated tensor. Unlike the sparse and low-rank settings directly on the frontal slices \citep{2020L2RM}, we characterize the inherited structures based on \blue{the low-rank, orthogonal decomposition (LROD). Similar to CPD}, the orthogonal decomposition can also reduce the model dimensionality, thereby reducing the computational complexity of the regression procedure. In addition, the encouraging sparsity on factor matrices produced by decomposition has been shown not only to yield asymptotically consistent estimators in high-dimensional settings \citep{2009JASAPCA}, but also to simplify visualization and interpretation of data analysis results \citep{2012HPCA}. \blue{It is worth mentioning that there is always an optimal LROD approximation \citep{orthogonalCP} for a given tensor, while the rank-$k$ approximation problem for fitting a CPD tensor is ill-posed for many rank values of $k$ in general \citep{CPDNOS-2008}. This makes the LROD tensor regression procedure more stable than those CPD based methods.}

To investigate the selection and estimation properties, we develop a sparse and low-rank tensor regression procedure by formulating it as a gSLOPE penalized LROD tensor optimization problem, namely TgSLOPE. Then, to measure its feature selection performance in finite samples, we define the notion of the tensor group FDR (TgFDR). Under the column-orthogonality assumption on the predictor matrix, TgSLOPE is shown to control TgFDR at any given level $0<q<1$ with appropriate choice of the regularization parameters. Moreover, the tensor estimator produced by TgSLOPE provably achieves the asymptotically minimax rate. Overall, our proposed TgFDRT controlling method can also be minimax optimal with respect to the estimation risk.

To well resolve our proposed TgSLOPE model with Stiefel manifold constraints, we constructively reformulate this manifold optimization problem as a difference-of-convex (DC) program whose objective function is shown to be level-coercive. This allows us to adopt some globally convergent DC-type algorithm in which the decision variables are updated by the proximal operator of the gSLOPE penalty. Simulations on synthetic data verify the TgFDR control, and test the effects of sparsity, model size and LROD rank on TgSLOPE performance respectively. In addition, numerical results on both synthetic data and a real human brain connection data confirm the superiority of our proposed TgSLOPE procedure by comparing it with several state-of-the-art approaches.

\subsection{Notation and preliminaries}
\begin{figure}
\begin{center}
           \subfigure[Horizontal slices $\bm{A}_{j_1::}$]{
          \includegraphics[height=1.1in,width=1.4in]{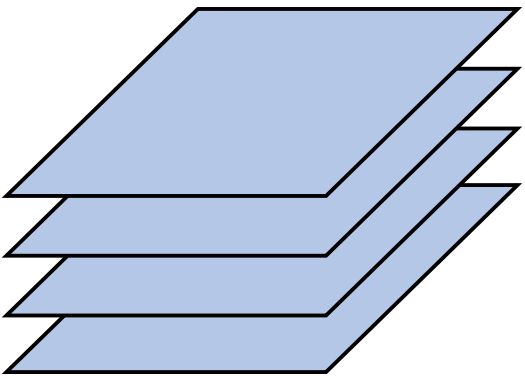}}
          \quad\quad
          \subfigure[Lateral slices $\bm{A}_{:j_2:}$]{
          \includegraphics[height=1.05in,width=1.3in]{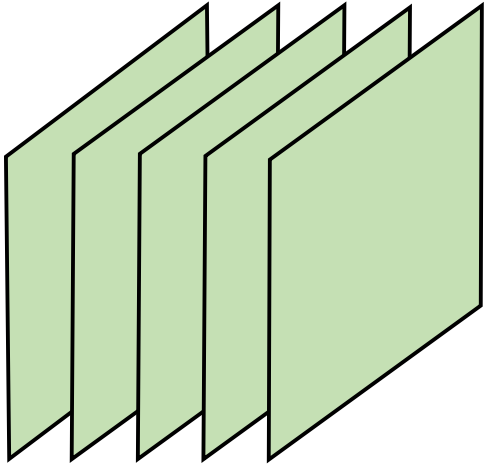}}
          \quad\quad
          \subfigure[Frontal slices $\bm{A}_{::j_3}$]{
          \includegraphics[height=1.15in,width=1.4in]{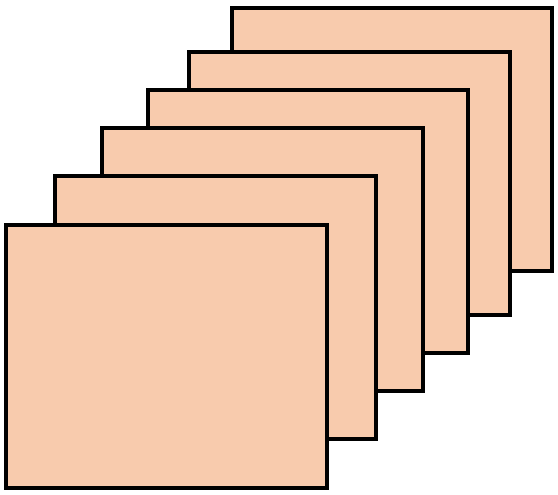}}
           \caption{\small{Slices of a order-3 tensor.}}\label{fig1}
\end{center}
\end{figure}

Throughout the article, we denote scalars, vectors, matrices and tensors by lowercase letters (e.g., $a, b$), boldface lowercase letters (e.g., $\bm{a}, \bm{b}$), boldface uppercase letters (e.g., $\bm{A}, \bm{B}$), and calligraphic letters (e.g., $\A, \B$), respectively. In addition, zero scalars, vectors, matrices and tensors are respectively denoted by $0, \bm{0}, \bm{O}$ and $\O$. For a vector $\bm{a}\in\mathbb{R}^n$, we denote $\bm{a}_{[1:s]} = (a_1,\dots,a_s)^\top$ with $s\leq n$. For a matrix $\bm{A}\in\mathbb{R}^{m\times n}$, the $i$th row and the $j$th column are denoted by $\bm{a}_{i:}$ and $\bm{a}_{:j}$, respectively, then write $\bm{A}_{[:,1:s]}=(\bm{a}_{:1},\dots,\bm{a}_{:s})$ with $s\leq n$. Figure \ref{fig1} shows the horizontal, lateral, and frontal slices of the tensor $\A\in\mathbb{R}^{p_1\times p_2\times p_3}$, denoted by $\bm{A}_{j_1::}\in\mathbb{R}^{p_2\times p_3}, \bm{A}_{:j_2:}\in\mathbb{R}^{p_1\times p_3}$, and $\bm{A}_{::j_3}\in\mathbb{R}^{p_1\times p_2}$, respectively. For convenience, we simply denote the $i$th row of $\bm{A}$ as $\bm{a}_i$ and the $j_3$th frontal slice of $\A$ as $\bm{A}_{j_3}$.

Given vectors $\bm{a}\in\mathbb{R}^m$ and $\bm{b}\in\mathbb{R}^n$, denote the outer product $\bm{a}\circ \bm{b}=\bm{a}\bm{b}^\top\in\mathbb{R}^{m\times n}$ and Kronecker product $\bm{a}\otimes\bm{b}=(a_1\bm{b}^\top,\dots,a_m\bm{b}^\top)^\top\in\mathbb{R}^{mn}$. For matrices $\bm{A}\in\mathbb{R}^{m\times q}$ and $\bm{B}\in\mathbb{R}^{n\times q}$, the Khatri-Rao product is defined as $\bm{A}\odot\bm{B}=(\bm{a}_{:1}\otimes\bm{b}_{:1},\dots,\bm{a}_{:q}\otimes\bm{b}_{:q})\in\mathbb{R}^{mn\times q}$; if $m=n$, the Hadamard product is defined as $\bm{A}*\bm{B}=[a_{ij}b_{ij}]\in\mathbb{R}^{m\times q}$, and the inner product is defined as $\langle\bm{A},\bm{B}\rangle=\sum_{i,j}a_{ij}b_{ij}$. For an order-3 tensor $\A\in\mathbb{R}^{p_1\times p_2\times p_3}$, the mode-$3$ product with a matrix $\bm{X}\in\mathbb{R}^{n\times p_3}$ is denoted by $\A\times_3\bm{X}\in \mathbb{R}^{p_1\times p_2\times n}$ with elements $(\A\times_3\bm{X})_{i,j,l}=\sum_{k=1}^{p_3}\A_{{i,j,k}}\bm{X}_{l,k}$. We also denote the matricization operator as $\M_3(\cdot)$, which unfolds the tensor $\A$ along the third mode into the matrix $\M_3(\A)\in\mathbb{R}^{p_3\times p_1p_2}$. Specifically, $(\M_3(\A))_{k,l}=\A_{i,j,k}$ with $l=1+(i-1)+(j-1)p_1$. Then the inverse of mode-$3$ unfolding can be denoted as $\M_3^{-1}(\cdot)$. The tensor $\A$ is called rank-one if it can be written as the outer product of vectors $\A = \bm{u} \circ \bm{v} \circ \bm{w}$. More about tensor operations can be found in \cite{2009DecomReview}.

Next we introduce some norms. For a vector $\bm{a}\in \mathbb{R}^{n}$, denote its $\ell_{2}$-norm as $\|\bm{a}\|=\sqrt{\sum_i a_i^2}$ and $\ell_{0}$-norm as $\|\bm{a}\|_{0} = {\rm \sharp}\{i: |a_i|\neq 0\}$. For a matrix $\bm{A}\in\mathbb{R}^{m\times n}$, denote its Frobenius norm as $\|\bm{A}\|_F=\sqrt{\sum_{i,j} a_{ij}^2}$ and trace as $\tr(\bm{A})=\sum_{i=j} a_{ij}$. Write singular values of $\bm{A}$ as $\sigma_1(\bm{A})\geq \dots \geq \sigma_p(\bm{A})$ with $p=\min\{m,n\}$. Then the nuclear norm is $\|\bm{A}\|_*=\sum_k\sigma_k(\bm{A})$ and the spectral norm is $\|\bm{A}\|_2=\sigma_1(\bm{A})$. We also introduce the notation $\|\bm{A}\|_r = (\|\bm{a}_{1}\|, \dots,\|\bm{a}_{m}\|)^\top$ for any matrix $\bm{A}\in\mathbb{R}^{m\times n}$ and $\|\A\|_f = (\|\bm{A}_1\|_F, \dots,\|\bm{A}_{p_3}\|_F)^\top$ for any tensor $\A\in\mathbb{R}^{p_1\times p_2\times p_3}$.

Moreover, we say that a random variable $X$ has a chi distribution with $n$ degrees of freedom, written by $X\sim \chi_n$, if it is the square root of a chi-squared random variable, i.e., $X^2\sim \chi^2_n$. For a proper closed convex function $f:\mathbb{R}^{m\times n}\rightarrow (-\infty,\infty]$, the subdifferential of $f$ at any given $\bm{X}\in dom(f)$, says $\partial f(\bm{X})$, is defined by
$$\partial f(\bm{X})=\{\bm{G}\in\mathbb{R}^{m\times n}: f(\bm{Y})\geq f(\bm{X})+\langle\bm{G}, \bm{Y}-\bm{X}\rangle \ for \ all \ \bm{Y}\in\mathbb{R}^{m\times n}\},$$ where each matrix $\bm{G}\in \partial f(\bm{X})$ is called a subgradient of $f$ at $\bm{X}$. The proximal operator of $f$ is given by
$$\prox_{f}(\bm{X})=\argmin\limits_{\bm{Z}}\Big\{f(\bm{Z})+\frac{1}{2}\|\bm{Z}-\bm{X}\|_F^2\Big\}, ~~\forall \bm{X}\in\mathbb{R}^{m\times n}.$$
Given positive scalars $a$ and $b$, denote $a\sim b$ if $a/b\to 1$, and $a \asymp b$ if there exist uniform constants $c,C>0$ such that $ca\leq b\leq Ca$.

The remainder of this article is organized as follows. In Section \ref{sec:2}, we introduce a frontal slice sparse and low-rank tensor regression procedure, which optimizes the gSLOPE penalized LROD tensor optimization problem, \blue{and then the identifiability of parameters is analyzed.} Section \ref{sec:3} makes statistical theory analysis for TgSLOPE procedure, including the TgFDR control at the prespecified level and the asymptotically minimax convergence with respect to the $\ell_2$-loss. An efficient and  globally convergent pDCAe algorithm is proposed in Section \ref{sec:4}. Section \ref{sec:5} reports some numerical studies to verify the performance of our proposed TgSLOPE approach. Concluding remarks are drawn in Section \ref{sec:6}. All proofs are deferred to the Appendix.

\section{Model} \label{sec:2}
In this section, we propose a feature selection and tensor estimate method for the tensor regression model (\ref{tm-1}), which can also be rewritten as
\begin{eqnarray}\label{tm-2}
\Y=\B^*\times_3\bm{X}+\E,
\end{eqnarray}
where $\Y, \E\in\mathbb{R}^{p_1\times p_2\times n}, \bm{X}\in\mathbb{R}^{n\times p}$. As discussed above, it is crucial to introduce some sparse and low-rank structures in order to facilitate estimation of the ultrahigh dimensional unknown parameters in finite samples. Thus, for the efficient feature selection, we consider the frontal slice sparsity of the coefficient tensor $\B^*$ based on the low-rank decomposition in the tensor regression framework (\ref{tm-2}).

\subsection{TgSLOPE estimate}
\begin{figure}
\begin{center}
          \includegraphics[width=5in]{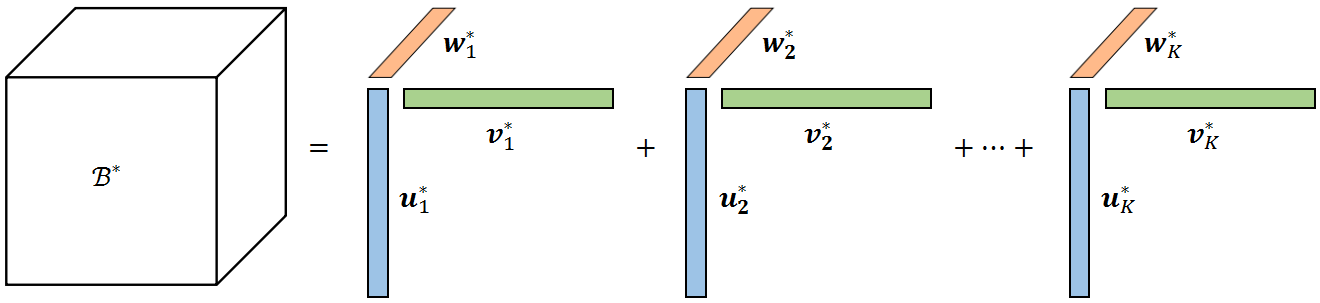}
           \caption{\small{Illustration of the rank-$K$ CPD of the tensor $\B^*$.}}\label{fig2}
\end{center}
\end{figure}

Assume that the coefficient tensor is LROD, that is, the tensor $\B^*$ admits a rank-$K$ CPD with column-orthogonal factor matrices. Specifically, the rank-$K$ CPD models a tensor as a sum of $K$ rank-one tensors \citep{2009DecomReview}, i.e., $\B^*=\sum_{k=1}^{K}\bm{u}_{k}^*\circ\bm{v}_{k}^*\circ\bm{w}_{k}^*$,
where the factor vectors $\bm{u}_{k}^*\in\mathbb{R}^{p_1}, \bm{v}_{k}^*\in\mathbb{R}^{p_2}, \bm{w}_{k}^*\in\mathbb{R}^{p}, k=1,\dots,K$. Figure \ref{fig2} illustrates the rank-$K$ CPD.
Denote the factor matrices composed of the vectors from the rank-one components by $\bm{U}^*=(\bm{u}_{1}^*, \bm{u}_{2}^*,\dots,\bm{u}_{K}^*)\in\mathbb{R}^{p_1\times K}$, $\bm{V}^*\in\mathbb{R}^{p_2\times K}$ and $\bm{W}^*\in\mathbb{R}^{p\times K}$ respectively. The CPD can be rewritten as $\B^*=[\![\bm{U}^*,\bm{V}^*,\bm{W}^*]\!]$. 
The LROD assumes that the factor matrices $\bm{U}^*$ and $\bm{V}^*$ produced by CPD are column-orthogonal. We define the rank of LROD tensor as the CP rank $K$. Moreover, we consider the row sparsity settings on the mode-3 factor matrix $\bm{W}^*$, which can translate to
the sparsity of the frontal slices for the tensor $\B^*$.

To identify significant features and estimate the coefficient tensor $\B^*$, we develop a penalized sparse and low-rank tensor regression method which optimizes the following gSLOPE penalized LROD tensor optimization problem (TgSLOPE)
\begin{eqnarray}\label{cp1-gslope}
\begin{split}
&\min\limits_{\bm{U},\bm{V}, \bm{W}}\quad \frac{1}{2}\big\|\Y-[\![\bm{U},\bm{V},\bm{W}]\!]\times_3\bm{X}\big\|_F^2+P_{\bm{\lambda}}\big(\|\bm{W}\|_r\big)\\
&\quad \ \ {\rm s.t.} \quad\quad \bm{U}^\top\bm{U}=\bm{V}^\top\bm{V}=\bm{I}_{K}, 
\end{split}
\end{eqnarray}
where $P_{\bm{\lambda}}(\bm{x})=\sum_{j=1}^p\lambda_j|x|_{(j)}$ is the SLOPE penalty with the regularization parameter vector $\bm{\lambda}=(\lambda_1,\dots,\lambda_p)^\top$ satisfying
\begin{equation*}\label{lambda_component} \lambda_1\geq\lambda_2\geq\cdots\geq\lambda_{p}\geq0  ~\mbox{and}~\lambda_1>0, \end{equation*}
and $|x|_{(j)}$ the $j$th largest component of $\bm{x}$ in magnitude.

According to the property of the matrix Khatri-Rao product that for any $\bm{A}\in\mathbb{R}^{m\times k}$ and $\bm{B}\in\mathbb{R}^{n\times k}$, $(\bm{A}\odot \bm{B})^\top(\bm{A}\odot \bm{B})=(\bm{A}^\top\bm{A})*(\bm{B}^\top\bm{B})$, we have the following lemma.
\begin{lemma}\label{lemma-1}
For any two column-orthogonal matrices $\bm{A}\in\mathbb{R}^{m\times k}$ and $\bm{B}\in\mathbb{R}^{n\times k}$, we have $\bm{C}=\bm{A}\odot \bm{B}\in \mathbb{R}^{mn\times k}$ is column-orthogonal, that is, $\bm{C}^\top\bm{C}=\bm{I}_k$.
\end{lemma}

Following the matricized form of a tensor, we have $\M_3(\B^*) = \bm{W}^*(\bm{V}^*\odot\bm{U}^*)^\top$. Denote 
$\bm{H}^*=\bm{V}^*\odot\bm{U}^*\in\mathbb{R}^{p_1p_2\times K}$. Without loss of generality, we assume that $K\leq \min\{p,p_1p_2\}$. It is known from Lemma \ref{lemma-1} that $\bm{H}^{*\top}\bm{H}^*=\bm{I}_K$. Thus, the TgSLOPE problem (\ref{cp1-gslope}) can be simplified as
\begin{eqnarray}\label{cp2-gslope}
\begin{split}
&\min_{\bm{W},\bm{H}}\quad L(\bm{W}, \bm{H})+P_{\bm{\lambda}}\big(\|\bm{W}\|_r\big)\\
&\ {\rm s.t.} \quad\ \bm{H}^\top\bm{H}=\bm{I}_{K},
\end{split}
\end{eqnarray}
where the loss $L(\bm{W}, \bm{H})=\frac{1}{2}\big\|\M_3(\Y)-\bm{X}\bm{W}\bm{H}^\top\big\|_F^2$.
The estimator $(\hat{\bm{W}}, \hat{\bm{H}})$ produced by (\ref{cp2-gslope}) and the tensor estimator $\hat{\B}$ are linked via $\M_3(\hat{\B})= \hat{\bm{W}}\hat{\bm{H}}^\top$.

\blue{Notably, the LROD tensor decomposition with the additional column-orthogonality on $U$ and $V$ turns to be a special CPD. While the best CP low-rank approximation of a tensor may not exist and the exact CP rank of a tensor is hard to compute, the imposed column-orthogonality on $\bm{U}$ and $\bm{V}$ can make sure that there is always an optimal rank-$k$ approximation for a given tensor \citep{orthogonalCP}. The resulting approximation tensor that admits such a special CPD has been coined as an orthogonally decomposable tensor, which has been studied in the communities of tensor analysis and tensor learning, see, e.g., the perturbation bounds and the applications to the unsupervised learning scenario of tensor SVD and the supervised task of tensor regression by \cite{LROD-Yuan}, and the linear convergence of an alternating polar decomposition method for low rank orthogonal tensor approximations by \cite{Hu-MP2023}. Additionally, \cite{2022Orthogonal} have illustrated that the orthogonally decomposable tensor model has the potential to perform better than CPD model in terms of predictive performance and model interpretation in numerical experiments. Technically, the column-orthogonal constraints on factor matrices can help to fix the indeterminacy and the non-uniqueness of CPD, thereby improving identifiability results of CPD methods (see, e.g., \cite{2013CP}).}

\blue{\subsection{Identifiability}
The identifiability of parameters is analyzed in this subsection. Considering CPD of a tensor $\B=[\![\bm{U},\bm{V},\bm{W}]\!]$, \cite{2013CP} have stated that the parameters in the tensor model is not identifiable due to two complications. One is the indeterminacy coming from scaling and permutation, and the second is possibly non-uniqueness of decomposition. For LROD tensor, the scaling indeterminacy can be avoided automatically due to the column-orthogonal $\bm{U}$ and $\bm{V}$. To fix the permutation indeterminacy, we adopt the convention that the first row of the factor matrix $\bm{U}$ is assumed to be arranged in descending order. Moreover, it follows from the Kruskal's condition of uniqueness of CPD (Theorem 4b, \cite{UCP-1977}) that the LROD is unique up to permutation if ${\rm k}(\bm{W})\geq 2$, where ${\rm k}(\bm{W})$ is the k-rank of $\bm{W}$ defined by the largest integer $k$ such that every subset of $k$ columns of $\bm{W}$ is linearly independent. Under the assumptions of determinacy and uniqueness of LROD, the imposed frontal slice sparsity on coefficient tensor $\B$ can be equivalently translated to the row sparsity settings on the mode-3 factor matrix $\bm{W}$, which yields the feasible set of problem (\ref{cp2-gslope}) being
\begin{align}\label{sp-lrtr}
\mathbb{B}_T=\big\{\B=[\![\bm{U},\bm{V},\bm{W}]\!]: \bm{U}^\top\bm{U}=\bm{V}^\top\bm{V}=\bm{I}_{K}, \big\|\|\bm{W}\|_r\big\|_0\leq s\big\},
\end{align}
where $s>0$ is a prescribed parameter that controls the frontal slice sparsity.

We give conditions of global identifiability result in the following proposition. The proof is presented in Appendix A.
\begin{proposition}\label{prop-identif}
Consider the tensor regression model (\ref{tm-2}) with the coefficient tensor $\B^*= [\![\bm{U}^*,\bm{V}^*,\bm{W}^*]\!]\in\mathbb{B}_T$. Then $\B^*$ is globally identifiable up to permutation if ${\rm k}(\bm{W}^*) \geq 2$ and $\big\|\|\B^*\|_f\big\|_0 \leq {\rm k}(\bm{X})/2$.
\end{proposition}
\begin{remark}
For sparse parameter models, \cite{Donoho-Identifiable-2003} have discussed the conditions for the uniqueness of sparse coefficient vector in classical linear regression frameworks. The identifiability result in Proposition \ref{prop-identif} can be regarded as an extension of Corollary 1 \citep{Donoho-Identifiable-2003} to the sparse and decomposable tensor regression models. Provided that the LROD tensor decomposition is unique up to permutation, the full coefficient tensor $\B^*$ with the frontal slice sparsity is globally identifiable and thus the LROD of $\B^*$ is identifiable.
\end{remark}

}
\section{Statistical Results \label{sec:3}}
This section is devoted to the TgFDR control and the estimate accuracy for our proposed TgSLOPE procedure.
\subsection{TgFDR Control}
FDR is a commonly used error rate that counts the expected proportion of errors among the rejected hypotheses in multiple testing. In this subsection, we develop the classic FDR notion to the setting of tensor regression and show that it can be well controlled by our proposed TgSLOPE.
\begin{definition}\label{def-tfdr}
Consider the tensor regression model (\ref{tm-2}) and let $(\hat{\bm{W}}, \hat{\bm{H}})$ be an estimator given by the optimization problem (\ref{cp2-gslope}). We define the tensor group false discovery rate (\TgFDR) for TgSLOPE as
\begin{align}\label{fdr-power}
\TgFDR = \mathbb{E}\bigg[\frac{V}{\max\{R,1\}}\bigg],
\end{align}
where $V, R$ are defined as follows
\begin{align*}
&V = \sharp\{j\in[p]: \bm{W}^*_j=\bm{0}, \hat{\bm{W}}_j\neq\bm{0}\}, \quad R = \sharp\{j\in[p]: \hat{\bm{W}}_j\neq\bm{0}\}
\end{align*}
with $\bm{W}_j^*$ and $\hat{\bm{W}}_j$ the $j$th rows of $\bm{W}^*$ and $\hat{\bm{W}}$, respectively.
\end{definition}
\blue{Define the regularization parameters of the TgSLOPE procedure as
\begin{align}\label{tuning}
\lambda_j=\sigma F_{\chi_{K}}^{-1}\big(1-q\cdot j/p\big), \ j\in[p],
\end{align}
where $0<q<1$, $F_{\chi_{K}}^{-1}(\alpha)$ is the $\alpha$th quantile of the $\chi$ distribution with K degrees of freedom. We give an upper bound of TgFDR in the following theorem. The technical proof is deferred to Appendix B.1.}
\begin{theorem}\label{th-tfdr}
Consider the tensor regression model (\ref{tm-2}) with the predictor matrix $\bm{X}$ satisfying $\bm{X}^\top\bm{X}=\bm{I}_p$. Then, for any solution $(\hat{\bm{W}}, \hat{\bm{H}})$ given by the TgSLOPE problem (\ref{cp2-gslope}) \blue{with the regularization parameters in (\ref{tuning})}, $\TgFDR$ obeys
\begin{align*}
\TgFDR = \mathbb{E}\bigg[\frac{V}{\max\{R,1\}}\bigg]\leq q\cdot \frac{p-s}{p}
\end{align*}
with $s$ the number of nonzero rows of $\bm{W}^*$.
\end{theorem}

\begin{figure}
\begin{center}
          \includegraphics[width=6in]{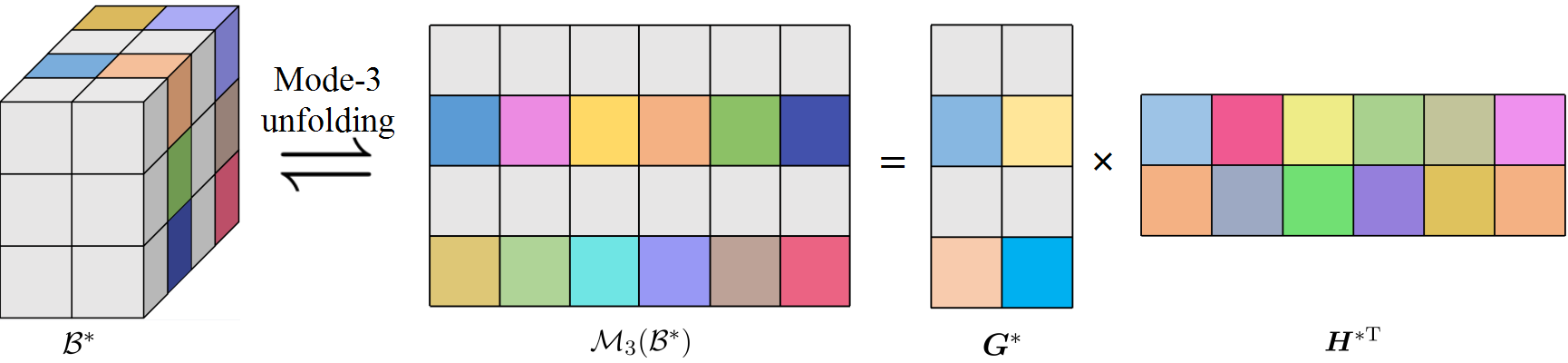}
           \caption{\small{Illustration of the sparsity equivalence relation between factor matrix $\bm{W}^*$ and tensor $\B^*$. Here the elements of matrices and tensors are zero (gray) and nonzero (color).}}\label{figadd}
\end{center}
\end{figure}

\begin{remark}
Under the guarantee of the uniqueness for LROD, our row sparsity settings on the matrix $\bm{W}^*$ can be equivalently interpreted as the sparsity in the frontal slices of the coefficient tensor $\B^*$. Figure \ref{figadd} illustrates this equivalence relation. Therefore, the definition of TgFDR in (\ref{fdr-power}) can be redefined based on the following random variables
\begin{align*}
TV = \sharp\{j\in[p]: \bm{B}_j^*=\bm{O}, \hat{\bm{B}}_j\neq\bm{O}\}, \quad TR = \sharp\{j\in[p]: \hat{\bm{B}}_j\neq\bm{O}\}
\end{align*}
with $\bm{B}_j^*$ and $\hat{\bm{B}}_j$ the $j$th frontal slices of $\B^*$ and $\hat{\B}$, respectively. This indicates that the estimator $\hat{\B}=\M_3^{-1}(\hat{\bm{W}}\hat{\bm{H}}^\top)$ given by (\ref{cp2-gslope}) provably controls TgFDR at any prespecified level $q\in(0,1)$.
\end{remark}
 \subsection{Estimate Accuracy}
This subsection aims to show that TgSLOPE enjoys minimax optimal estimation property under the assumption that the ground truth $\B^*$ is bounded, that is, $\max\{\Vert\bm{B}_j^*\Vert_F, j\in[p]\}< \infty$. We measure the deviation of an estimator from $\B^*$ in the following theorem. The technical proof is given in Appendix B.2.
 \begin{theorem}\label{thm-2}
 Consider the TgSLOPE optimization problem (\ref{cp2-gslope}) with the predictor matrix $\bm{X}$ satisfying $\bm{X}^\top\bm{X}=\bm{I}_p$ and the regularization parameters 
 in (\ref{tuning}). Then, under $p\to \infty$ and $s/p\to 0$, the TgSLOPE estimator \blue{
 enjoys the asymptotically minimax rate over $\mathbb{B}_{s}=\{\B: \big\|\|\B\|_f\big\|_{0}\leq s\}$, that is,
 \begin{align}\label{sp-1}
 \sup_{\B^*\in \mathbb{B}_{s}}\mathbb{E}\|\hat{\B}-\B^*\|_F^2  \sim
\inf_{\tilde{\B}}\sup_{\B^*\in \mathbb{B}_{s}}\mathbb{E}\Vert \tilde{\B}-\B^* \Vert_F^2,
 \end{align}
 }
where the infimum is taken over all estimators $\tilde{\B}$ based on date set $\D=\{(\bm{x}_i,\bm{Y}_i)\}_{i=1}^n$.
 \end{theorem}
\blue{
\begin{remark}
The asymptotically minimax result in Theorem \ref{thm-2} is an extension of Theorem 1.1 in \cite{2016slope} (for the SLOPE estimate) and Theorem 2.2 in \cite{2019gslope} (for the group SLOPE estimate) under the classical linear regression frameworks. In addition, we note that the choice of the regularization parameters in Theorem \ref{thm-2} does not depend on the sparsity level. In such a sense, TgSLOPE is adaptive to a range of sparsity in achieving the minimax optimality. It would also be interesting to work on estimation bounds of minimax rates in a non-asymptotic manner,  such as that in \cite{2019Convex} and \cite{2020Cubic}. We leave this as one of our research topics in the future.
\end{remark}
}

\section{Proximal DCA with Extrapolation \label{sec:4}}
In this section, we propose a proximal difference-of-convex algorithm with extrapolation (pDCAe) to solve the TgSLOPE optimization problem and establish its global convergence.

\subsection{DC Reformulation}
Given an optimal solution $(\hat{\bm{W}}, \hat{\bm{H}})$ of the problem (\ref{cp2-gslope}), it is known from the orthogonal Procrustes problem \citep{Procrustes} that $\hat{\bm{H}}$ can be determined by $\hat{\bm{W}}$ in the way that $\hat{\bm{H}}=\hat{\bm{U}}_{[:,1:K]}\hat{\bm{V}}^\top$, where $\hat{\bm{U}}, \hat{\bm{V}}$ are obtained from the SVD $\M_3(\Y)^\top\bm{X}\hat{\bm{W}}= \hat{\bm{U}}\bm{D}\hat{\bm{V}}^\top$ with $\hat{\bm{U}}\in\mathbb{R}^{p_1p_2 \times p_1p_2}$, $\bm{D}\in\mathbb{R}^{p_1p_2 \times K}$ and $\hat{\bm{V}}\in\mathbb{R}^{K \times K}$. Plugging $(\hat{\bm{W}}, \hat{\bm{H}})$ into the objective function of problem (\ref{cp2-gslope}), we obtain that
\begin{align*}
L(\hat{\bm{W}}, \hat{\bm{H}})
&=\frac{1}{2}\big\|\M_3(\Y)\big\|_F^2
+\frac{1}{2}\big\|\bm{X}\hat{\bm{W}}\hat{\bm{H}}^\top\big\|_F^2-\tr(\M_3(\Y)^\top\bm{X}\hat{\bm{W}}\hat{\bm{H}}^\top)\\
& = \frac{1}{2}\big\|\M_3(\Y)\big\|_F^2
+\frac{1}{2}\big\|\bm{X}\hat{\bm{W}}\big\|_F^2-\tr(\bm{D}\hat{\bm{U}}_{[:,1:K]}^\top\hat{\bm{U}})\\
& =\frac{1}{2}\big\|\M_3(\Y)\big\|_F^2
+\frac{1}{2}\big\|\bm{X}\hat{\bm{W}}\big\|_F^2-\tr(\bm{D})\\
&= \frac{1}{2}\big\|\M_3(\Y)\big\|_F^2
+\frac{1}{2}\big\|\bm{X}\hat{\bm{W}}\big\|_F^2-\|\M_3(\Y)^\top\bm{X}\hat{\bm{W}}\|_{*}.
\end{align*}
Therefore, the TgSLOPE problem (\ref{cp2-gslope}) can be equivalently transformed into
\begin{align}\label{dcpg-1}
\left\{
\begin{array}{l}
\hat{\bm{W}}=\argmin \Big\{F(\bm{W})=\frac{1}{2}\|\bm{X}\bm{W}\|_F^2+P_{\bm{\lambda}}\big(\|\bm{W}\|_r\big)-\|\M_3(\Y)^\top\bm{X}\bm{W}\|_{*}\Big\},\vspace{2ex}\\
\hat{\bm{H}}=\hat{\bm{U}}_{[:,1:K]}\hat{\bm{V}}^\top, \mbox{where $\hat{\bm{U}}$ \mbox{and} $\hat{\bm{V}}$ are from the SVD}\ \M_3(\Y)^\top\bm{X}\hat{\bm{W}}= \hat{\bm{U}}\bm{D}\hat{\bm{V}}^\top.
\end{array}
\right.
\end{align}
Denote $F1(\bm{W})=\frac{1}{2}\|\bm{X}\bm{W}\|_F^2$ and $F2(\bm{W}) = \frac{L}{2}\|\bm{W}\|_F^2- F1(\bm{W})$. Since $F1$ is gradient Lipschitz continuous with modulus $L=\|\bm{X}^\top\bm{X}\|_2$, the convexity of $F2$ follows. Inspired by the DC-representable function from \cite{2018pdca}, we rewrite the optimization problem in (\ref{dcpg-1}) as
\begin{align}\label{pdca-2}
\min_{\bm{W}\in\mathbb{R}^{p\times K}} F(\bm{W})= \underbrace{\frac{L}{2}\|\bm{W}\|_F^2+P_{\bm{\lambda}}\big(\|\bm{W}\|_r\big)}_{C1(\bm{W})}-\underbrace{\bigg(\frac{L}{2}\|\bm{W}\|_F^2-\frac{1}{2}\|\bm{X}\bm{W}\|_F^2+N(\bm{W})\bigg)}_{C2(\bm{W})},
\end{align}
where $N(\bm{W})=\|\M_3(\Y)^\top\bm{X}\bm{W}\|_{*}$. It is easy to verify that both $C1$ and $C2$ are convex, leading to a DC program as in (\ref{pdca-2}).

By now, we constructively reformulate the manifold optimization problem (\ref{cp2-gslope}) as a DC program (\ref{pdca-2}). This allows us to solve the reformulation problem of TgSLOPE by a DC-type algorithm.

\subsection{Proximal DCA}
Considering the DC program (\ref{pdca-2}), we first give the subdifferential of $N(\bm{W})$ in the following lemma, which can be clearly derived from the subdifferential of the matrix nuclear norm \citep{1992nuclear} and the chain rule of Theorem 10.6 in \cite{VA}.
\begin{lemma}\label{lemma-2}
The subdifferential of $N(\bm{W})$ is given by
\begin{align}\label{sd-nunorm}
\partial N(\bm{W})=\{\bm{X}^\top\M_3(\Y)(\tilde{\bm{U}}_{[:,1:K]}\tilde{\bm{V}}^\top+\tilde{\bm{W}}): \tilde{\bm{U}}_{[:,1:K]}^\top\tilde{\bm{W}}=\bm{O}, \tilde{\bm{W}}\bar{\bm{V}}=\bm{O}, \|\tilde{\bm{W}}\|_2\leq 1\},
\end{align}
where $\tilde{\bm{U}}, \tilde{\bm{V}}$ are obtained by the SVD of $\M_3(\Y)^\top\bm{X}\bm{W}= \tilde{\bm{U}}\bm{D}\tilde{\bm{V}}^\top$.
\end{lemma}

The DCA iterative scheme for (\ref{pdca-2}) takes the following form
\begin{align}\label{pdca-3}
\begin{split}
\bm{W}^{(k+1)}&=\argmin_{\bm{W}}\bigg\{\frac{L}{2}\|\bm{W}\|_F^2-\Big\langle\bm{W}, (L\bm{I}-\bm{X}^\top\bm{X})\bm{W}^{(k)}+\bm{Q}_1^{(k)}\Big\rangle+P_{\bm{\lambda}}\big(\|\bm{W}\|_r\big)\bigg\}\\
&=\argmin_{\bm{W}}\bigg\{\frac{1}{2}\bigg\|\bm{W}-\bigg(\bm{W}^{(k)}-\frac{1}{L}\Big(\bm{X}^\top\bm{X}\bm{W}^{(k)}-\bm{Q}_1^{(k)}\Big)\bigg)\bigg\|_F^2
+P_{\bm{\lambda/L}}\big(\|\bm{W}\|_r\big)\bigg\},
\end{split}
\end{align}
where $\bm{Q}_1^{(k)}\in \partial N(\bm{W}^{(k)})$. Then the optimal solution of problem (\ref{pdca-3}) can be computed by the proximal operator of the penalty function $P_{\bm{\lambda}/L}$. Specifically,
\begin{align}\label{pdca-4}
\left\{
\begin{array}{l}
\bm{\eta}^{(k+1)}=\argmin_{\bm{\eta}}\Big\{\frac{1}{2}\big\|\|\bm{Q}^{(k)}\|_r-\bm{\eta}\big\|^2+P_{\bm{\lambda}/L}(\bm{\eta})\Big\},\vspace{1ex}\\
\bm{W}^{(k+1)}_j=\big(\prox_{P_{\bm{\lambda}/L}}(\bm{Q}^{(k)})\big)_j=\eta_j^{(k+1)}\frac{\bm{q}^{(k)}_j}{\|\bm{q}^{(k)}_j\|}, j\in[p],
\end{array}
\right.
\end{align}
where $\bm{Q}^{(k)}=\bm{W}^{(k)}-\big(\bm{X}^\top\bm{X}\bm{W}^{(k)}-\bm{Q}_1^{(k)}\big)/L$. The resulting DCA is called as the proximal DCA \citep{2018pdca,Chen-pdca}. Consequently, the $\bm{W}$-subproblem is reduced to identifying the general SLOPE proximal operator, which can be efficiently solved by FastProxSL1 \citep{2015slope}. Furthermore, in order to possibly accelerate the algorithm, we entertain extrapolation techniques \citep{Nesterov2013} in the proximal DCA. The algorithmic framework is summarized in Algorithm \ref{Alg2}.

\begin{algorithm}
\caption{pDCAe for solving TgSLOPE (\ref{cp2-gslope})\label{Alg2}}
 Initialize $\bm{W}^{(0)}\in \mathbb{R}^{p\times K}$, $\{\beta^{(k)}\} \subseteq [0,1)$ with $\sup_{k}\beta^{(k)}<1$, set $\bm{W}^{(-1)}= \bm{W}^{(0)}$, $k=0$;
\begin{algorithmic}[1]
\State Choose $\bm{Q}_1^{(k)}\in \partial N(\bm{W}^{(k)})$ by (\ref{sd-nunorm});

\State Compute $\bm{A}^{(k)} = \bm{W}^{(k)}+\beta^{(k)}(\bm{W}^{(k)}-\bm{W}^{(k-1)})$ and update
$$\bm{Q}^{(k)}=\bm{A}^{(k)}-\big(\bm{X}^\top\bm{X}\bm{A}^{(k)}-\bm{Q}_1^{(k)}\big)/L;$$

\State  Compute $\bm{W}^{(k+1)}$ by the two-step form in (\ref{pdca-4});

\State  Set $k=k+1$. If the stopping criterion is met, perform the SVD $\M_3(\Y)^\top\bm{X}\bm{W}^{(k)}=\bm{U}^{(k)}\bm{D}^{(k)}\bm{V}^{(k)\top}$ to get $$\bm{H}^{(k)}=\big(\bm{U}^{(k)}\big)_{[:,1:K]}\bm{V}^{(k)\top},$$ then stop and return $\B^{(k)}=\M_3^{-1}(\bm{W}^{(k)}\bm{H}^{(k)\top})$; otherwise, go to the step 1.
\end{algorithmic}
\end{algorithm}

Note that the pDCAe algorithm is reduced to the general proximal DCA if $\beta^{(k)}=0$ for all $k$. Moreover, Algorithm \ref{Alg2} can be coupled with many popular choices of extrapolation parameters $\{\beta^{(k)}\}$ including that used in accelerated proximal gradient (APG) for solving SLOPE \citep{2015slope}. In our numerical experiments section, we follow \cite{2015slope} and set $\theta^{(-1)}=1$,
$$\beta^{(k)} = \frac{\theta^{(k-1)}-1}{\theta^{(k)}}\quad {\rm with} \quad \theta^{(k)}=\frac{1+\sqrt{1+4(\theta^{(k-1)})^2}}{2}, \quad \forall k\geq 0.$$

\subsection{Global Convergence}
This subsection is dedicated to the global convergence of pDCAe. We start with the level-coercivity of the DC function $F(\bm{W})$.
\begin{lemma}\label{lemma-3}
Let the predictor matrix $\bm{X}$ in problem (\ref{cp2-gslope}) be full column rank. Then the DC function $F(\bm{W})$ in (\ref{pdca-2}) is level-coercive in the sense that $\liminf_{\|\bm{W}\|_F\to \infty}\frac{F(\bm{W})}{\|\bm{W}\|_F}>0.$
\end{lemma}

\blue{Based on the level-coercivity as obtained in Lemma \ref{lemma-3}, the desired global convergence of pDCAe for solving problem (\ref{cp2-gslope}) is stated as follows.}

\begin{theorem}\label{thm-3} \blue{Suppose that the prediction matrix $\bm{X}$ is full column rank and let} $\{\bm{W}^{(k)}\}$ be the sequence generated by Algorithm \ref{Alg2} for solving problem (\ref{cp2-gslope}). Then the following properties hold:
\begin{enumerate}[(a)]
\item  The sequence $\{\bm{W}^{(k)}\}$ is bounded;
\item  $\lim_{k\to \infty}\|\bm{W}^{(k+1)}-\bm{W}^{(k)}\|_F=0$;
\item  Every limit point $\bar{\bm{W}}$ of the sequence $\{\bm{W}^{(k)}\}$ is a stationary point of $F$ in the sense that
    $$\bm{O}\in \bm{X}^\top\bm{X}\bar{\bm{W}}+\partial P_{\bm{\lambda}}\big(\|\bar{\bm{W}}\|_r\big)-\partial \|\M_3(\Y)^\top\bm{X}\bar{\bm{W}}\|_{*}.$$
\end{enumerate}
\end{theorem}

The proofs of results in this subsection are provided in Appendix C.

\section{Numerical Experiments \label{sec:5}}
This section gives some experiments on synthetic data and a real human brain connection data. All numerical experiments are implemented in MATLAB (R2021a), running on a laptop with Intel Core i5-8265U CPU (1.60GHz) and 16 GB RAM.
\subsection{Comparative Methods}
We verify the performance of proposed TgSLOPE by comparing it with the following three approaches:
\begin{enumerate}[(a)]
\item TBMM: the block majorization minimization (BMM) algorithm proposed by \cite{2021SRRR} to solve the TgSLOPE problem (\ref{cp2-gslope}), whose corresponding iterative scheme is
\begin{align*}
\left\{
\begin{array}{l}
\bm{W}^{(k+1)}=\argmin_{\bm{W}}\big\{\frac{1}{2}\|\bm{W}-\bm{R}^{(k)}\|_F^2+P_{\bm{\lambda/L}}\big(\|\bm{W}\|_r\big)\big\},\vspace{1ex}\\
\bm{H}^{(k+1)}=\big(\bm{U}^{(k+1)}\big)_{[:,1:K]}\bm{V}^{(k+1)\top},
\end{array}
\right.
\end{align*}
where $\bm{R}^{(k)}=\bm{W}^{(k)}
-\big(\bm{X}^\top\bm{X}\bm{W}^{(k)}-\bm{X}^\top\M_3(\Y)\bm{H}^{(k)}\big)/L$, $\bm{U}^{(k+1)}$ and $\bm{V}^{(k+1)}$ are from the SVD $\M_3(\Y)^\top\bm{X}\bm{W}^{(k+1)}=\bm{U}^{(k+1)}\bm{D}^{(k+1)}\bm{V}^{(k+1)\top}$;

\item TgLASSO: the group LASSO penalized LROD tensor regression approach, in which the pDCAe algorithm in Section \ref{sec:4} is used to solve
\begin{align*}
\min_{\bm{W},\bm{H}} \bigg\{\frac{1}{2}\big\|\M_3(\Y)-\bm{X}\bm{W}\bm{H}^\top\big\|_F^2+\lambda\sum_{j=1}^p \|\bm{W}_j\|: \bm{H}^\top\bm{H}=\bm{I}_{K}\bigg\},
\end{align*}
where $\lambda>0$ is the tuning parameter;

\item TLRR: the LROD tensor regression (without sparse penalty), in which the pDCAe is applied to solve
\begin{align*}
\min_{\bm{W},\bm{H}} \bigg\{\frac{1}{2}\big\|\M_3(\Y)-\bm{X}\bm{W}\bm{H}^\top\big\|_F^2: \bm{H}^\top\bm{H}=\bm{I}_{K}\bigg\}.
\end{align*}
\end{enumerate}
\blue{In our numerical experiments, we simply adopt the stopping criterion proposed by \cite{Chen-pdca} for all approaches, i.e., $\frac{\|\bm{W}^{(k)}-\bm{W}^{(k+1)}\|_F}{\max\{\|\bm{W}^{(k)}\|_F,1\}}\leq \epsilon$ 
with some given tolerance $\epsilon>0$.}

\subsection{Synthetic Data}
In our numerical experiments on synthetic data, we adopt the measures including TgFDR defined in (\ref{fdr-power}) and the tensor power (TP) to evaluate selection performance of the estimator $\hat{\B}$ generated by a given method. Here TP is defined as
$${\rm TP} = \mathbb{E}(T)/s, \ \mbox{where} \  T = \sharp\{j\in[p]: \bm{B}_j^*\neq\bm{O}, \hat{\bm{B}}_j\neq\bm{O}\} \ \mbox{and} \ s = \sharp\{j\in[p]: \bm{B}_j^*\neq\bm{O}\}.$$
For estimation accuracy, we evaluate the performance of $\hat{\B}$ in terms of the relative group estimate error (RgEE) and the mean squared error (MSE) defined respectively as
$${\rm RgEE}=\frac{\big\Vert \|\hat{\B}\|_f-\|\B^*\|_f  \big\Vert^2}{\big\Vert\|\B^*\|_f  \big\Vert^2}, \quad {\rm MSE} =\big\| (\hat{\B}-\B^*)\times_3\bm{X} \big\|_F^2/np_1p_2.$$ Meanwhile, to evaluate the time efficiency of our proposed pDCAe algorithm, the CPU time (Time) is reported for each testing instance.

The ground truth is simulated via $\B^* = \M_3^{-1}(\bm{W}^*\bm{H}^{*\top})\in \mathbb{R}^{p_1\times p_2\times p}$, in which $\bm{W}^*\in\mathbb{R}^{p\times K}$ is generated in a similar manner as in \cite{2019gslope} with $s$ nonzero rows. Specifically, each nonzero row of $\bm{W}^*$ is generated from the uniform distribution $U[0.1,1.1]$ and then we scale it such that $\|\bm{w}_j^*\|=a\sqrt{K}$ with $a=\sqrt{4\ln(p)/(1-p^{-2/K})-K}$. \blue{The explanation of such a special simulation procedure is presented in Appendix D.} The column-orthogonal matrix $\bm{H}^*\in\mathbb{R}^{p_1p_2\times K}$ is simulated as the first $K$ left singular vectors of a $p_1p_2 \times p_1p_2$ matrix with i.i.d. standard normal distributed entries. The response tensor $\Y\in\mathbb{R}^{p_1\times p_2\times n}$ is simulated by $\Y=\B^*\times_3\bm{X}+\E$, where the entries of the noise tensor $\E$ are i.i.d. drawn from $N(0,1)$.

We first verify the TgFDR controlling performance of TgSLOPE under the following two situations of the predictor matrix $\bm{X}\in\mathbb{R}^{n\times p}$: \blue{(a) the orthogonal design, where $\bm{X}$ is generated from the orthogonalization for a $n \times p$ matrix with i.i.d. standard normal distributed entries; (b) the commonly used Gaussian random design, such as in \cite{2020L2RM}, where the predictor vectors $\bm{x}_i, i=1,\cdots, n$ are i.i.d. drawn from Gaussian distribution $\mathcal{N}(\bm{0},\bm{C})$, where $\bm{C}$ is a $p\times p$ covariance matrix with the element $\bm{C}_{j_1j_2}=0.5^{|j_1-j_2|}$ for any $j_1,j_2\in[p]$. The orthogonal situation works with the regularization parameter sequence defined in (\ref{tuning}). For Gaussian random design, we select regularization parameters of TgSLOPE according to Procedure 2 of \cite{2019gslope}.}

\begin{figure}
\begin{center}
           \subfigure[Orthogonal design]{
          \includegraphics[width=2.5in]{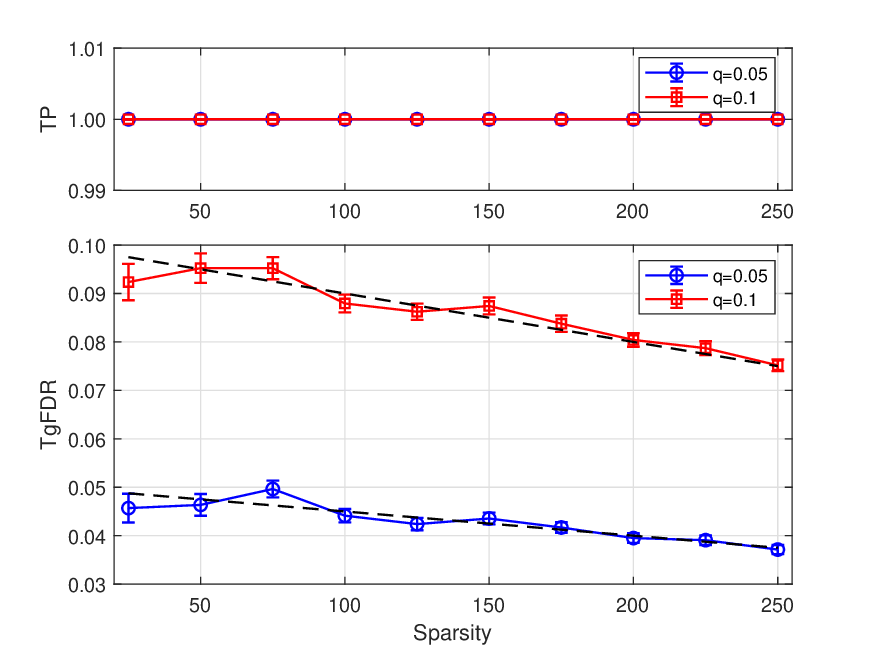}}
          \subfigure[\blue{Gaussian random design}]{
          \includegraphics[width=2.5in]{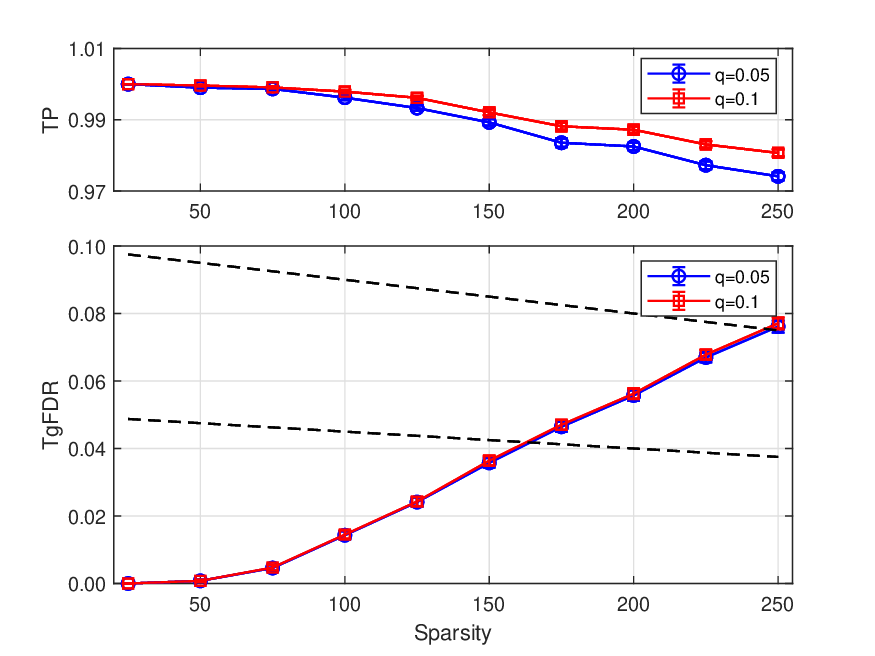}}
           \caption{\small{Estimated TgFDR and TP with \blue{$n=2000, p=1000$}, $p_1=p_2=10, K=20$ under the orthogonal and Gaussian situations of the matrix $\bm{X}$. Bars correspond to $\pm $SE (standard error), black dotted lines represent the `nominal' TgFDR level $q\cdot(p-s)/p$.}}\label{fig3}
\end{center}
\end{figure}

Set $n=2000, p=1000, p_1=p_2=10, K=20$ and sparsity $s=25:25:250$. We perform 100 independent replications for each sparsity and target TgFDR level ($q=0.05, 0.1$). As shown in Figure \ref{fig3} (a), TgSLOPE maintains a comparable TgFDR with the `nominal' level and reports the estimated tensor power TP =1 for all testing instances. \blue{For the Gaussian situation, TgSLOPE reports the relatively low TgFDR with the strong power, especially for the sparse cases as shown in Figure \ref{fig3} (b). The reported TgFDR of TgSLOPE does not exceed the `nominal' level of the case $q=0.1$ in all sparsity settings except for the case $s=250$.}

In the following numerical comparisons, we consider the \blue{Gaussian random situation in which the entries of $\bm{X}$ are i.i.d. drawn from $N(0,1/n)$}, and test the effect of sparsity $s$, model size $p$ and LROD rank $K$ respectively under the target TgFDR level $q=0.05$. For TgLASSO, the tuning parameter is selected by 5-fold cross-validation.

\begin{figure}
\begin{center}
           \subfigure[TgFDR]{
          \includegraphics[width=2.2in]{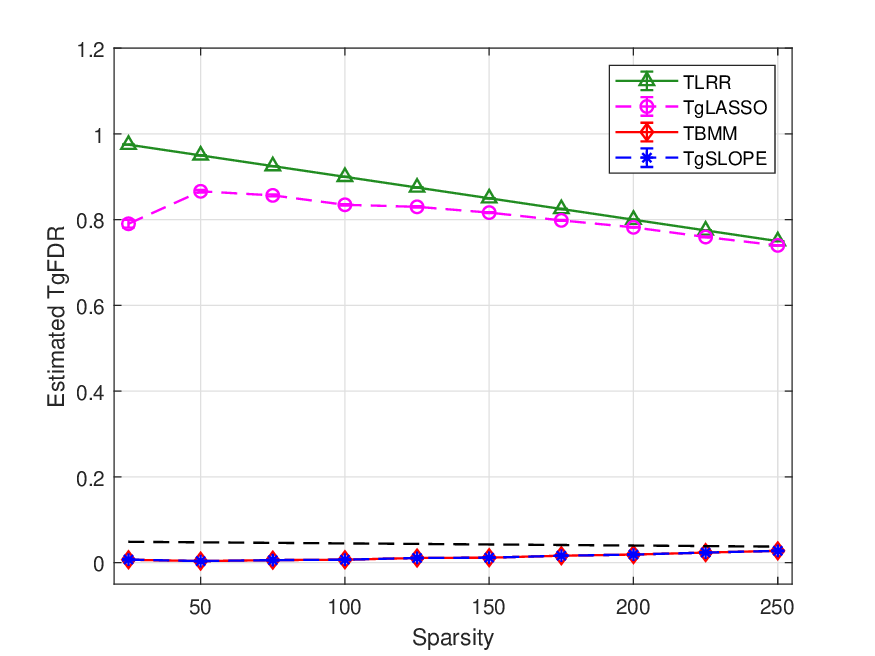}}
          \subfigure[RgEE]{
          \includegraphics[width=2.2in]{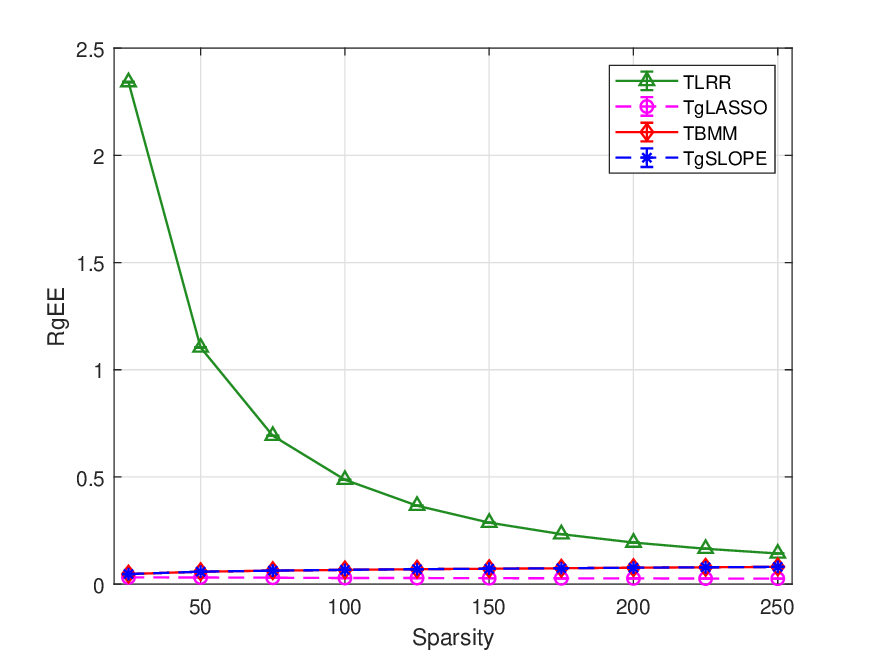}}
          \subfigure[MSE]{
          \includegraphics[width=2.2in]{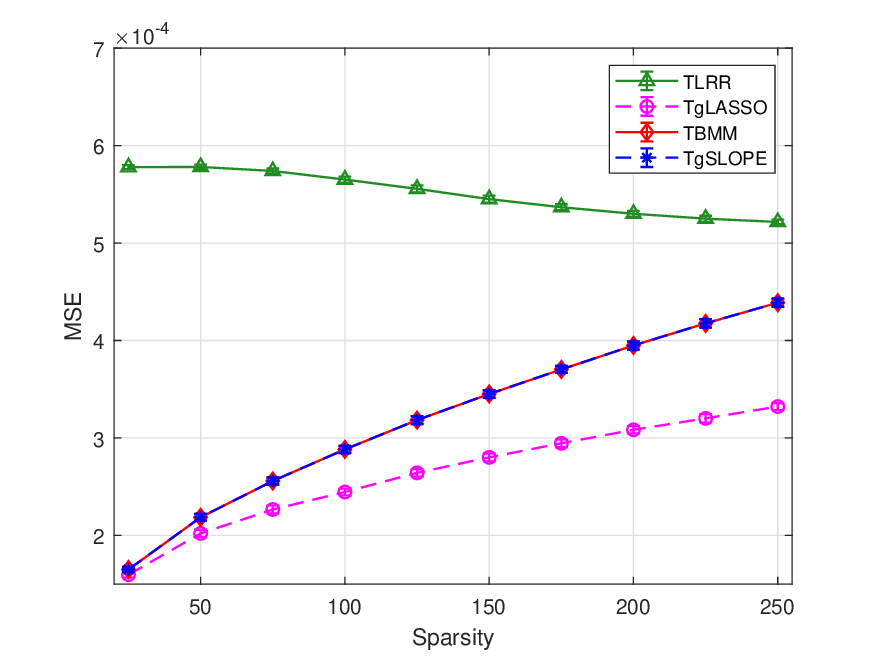}}
          \subfigure[Time]{
          \includegraphics[width=2.2in]{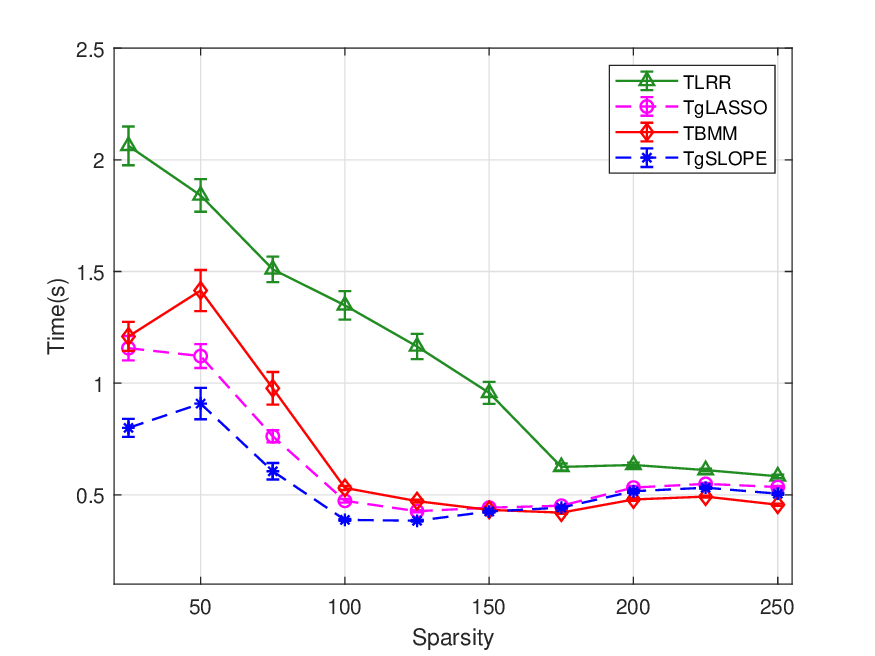}}
           \caption{Average results with various sparsity under $n=3000, p=1000, p_1=p_2=10, K=20$ and for the Gaussian random design situation of $\bm{X}$. Bars correspond to $\pm $SD (standard deviation). In (a), black dotted lines represent the `nominal' TgFDR level $q\cdot(p-s)/p$.}\label{fig4}
\end{center}
\end{figure}

Sparsity effect: set $n=3000, p=1000, p_1=p_2=10$ and $K=20$. We test the performance of four approaches under various sparsity with $s=25:25:250$. Simulation results report TP =1 for all competitors in each testing instance. In addition, average results based on 100 independent replications have been shown in Figure \ref{fig4}. (i) Figure \ref{fig4} (a) illustrates that TgSLOPE and TBMM have significant superiority in terms of TgFDR, with TgFDR below `nominal' for all testing sparsity. (ii) TgLASSO fails to control TgFDR, although it reports the relatively small RgEE and MSE as shown in Figure \ref{fig4} (b) and (c). (iii) Figure \ref{fig4} (d) shows that TgSLOPE reports a relatively short CUP time, especially in the cases of small sparsity. (iv) It is intuitive that those four methods tend to have comparable performance as the sparsity increases.

\begin{table}
\renewcommand{\arraystretch}{1.2}
\caption{Average results with different model size under $n=3000, p_1=p_2=10, K=20, s=0.02p$ for the Gaussian random design situation of $\bm{X}$. Standard deviations are presented in brackets.\label{tab-1}}
\setlength{\tabcolsep}{4mm}
\begin{center}
{\footnotesize \begin{tabular}{lllllll}
\hline
      Size   & Method    &TgFDR    &RgEE             &MSE                &Time(s)               \\ \hline

    $p=2000$  &TLRR       &0.980  &5.10e-0 (6.24e-2)   & 7.77e-4 (1.89e-6)   &3.325        \\
              &TgLASSO        &0.891  &3.44e-2 (3.97e-3)   & 1.92e-4 (4.45e-6)   &1.566    \\
              &TBMM   &0.006  &5.59e-2 (3.91e-3)   & 2.02e-4 (2.73e-6)   &2.242     \\
              &TgSLOPE &0.007  &5.59e-2 (3.91e-3)   & 2.02e-4 (2.76e-6)   &1.081    \\  \hline

     $p=4000$ &TLRR        &0.980   & 3.92e-0 (4.26e-2)   & 9.28e-4 (2.45e-6)   &10.163        \\
              &TgLASSO         &0.921   & 3.50e-2 (2.16e-3)   & 2.40e-4 (4.70e-6)   &3.727    \\
              &TBMM    &0.008   & 6.46e-2 (3.28e-3)   & 2.70e-4 (3.48e-6)   &6.104     \\
              &TgSLOPE  &0.008   & 6.45e-2 (3.27e-3)   & 2.70e-4 (3.49e-6)   &2.962    \\  \hline

     $p=6000$ &TLRR        &0.980  &1.69e-0 (1.71e-2)    &9.19e-4 (2.94e-6)  &15.146       \\
              &TgLASSO         &0.924   &3.66e-2 (2.17e-3)    &2.77e-4 (5.16e-6)  &6.179    \\
              &TBMM    &0.025   &6.91e-2 (3.78e-3)    &3.24e-4 (4.02e-6)  &9.108     \\
              &TgSLOPE  &0.014   &6.88e-2 (3.75e-3)    &3.24e-4 (4.01e-6)  &5.666    \\
\hline
 \end{tabular}}
\end{center}
\end{table}

Model size effect: in this example, we test the effect of model size $p$ comparing among four approaches. Set $n=3000, p_1=p_2=10, K=20$ and the sparsity $s=0.02p$ with model size $p=2000,4000,6000$. In each testing instance, all these competitors report TP =1. Moreover, Table \ref{tab-1} collects the average results based on 100 independent replications for each model size. (i) As suggested in Table \ref{tab-1}, TLRR gives the worst performance. For the other three approaches, all of the four evaluation metrics tend to become larger as the model size grows. (ii) Similar to the results shown in Figure \ref{fig4}, TgLASSO fails to control TgFDR, while TgSLOPE and TBMM can maintain the TgFDR below the `nominal' TgFDR level for all testing model sizes, with predictable sacrifice on the estimation accuracy. (iii) We can also see from Table \ref{tab-1} that TgSLOPE gives the smaller TgFDR than TBMM for the model size $p=6000$, and it reports the least CPU time for all cases.

\begin{figure}
\begin{center}
          \subfigure[TgFDR]{
          \includegraphics[width=2.2in]{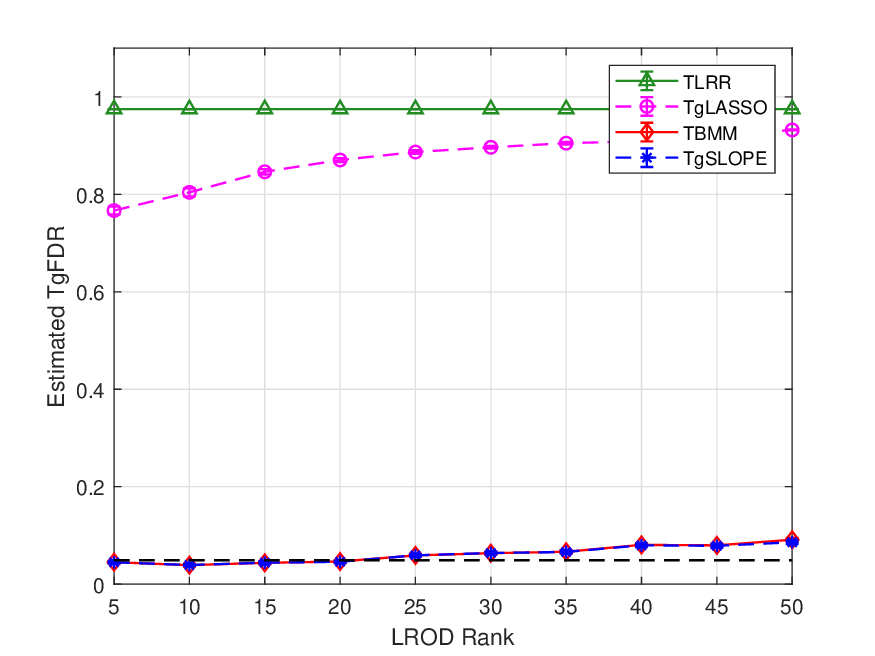}}
          \subfigure[TgFDR]{
          \includegraphics[width=2.2in]{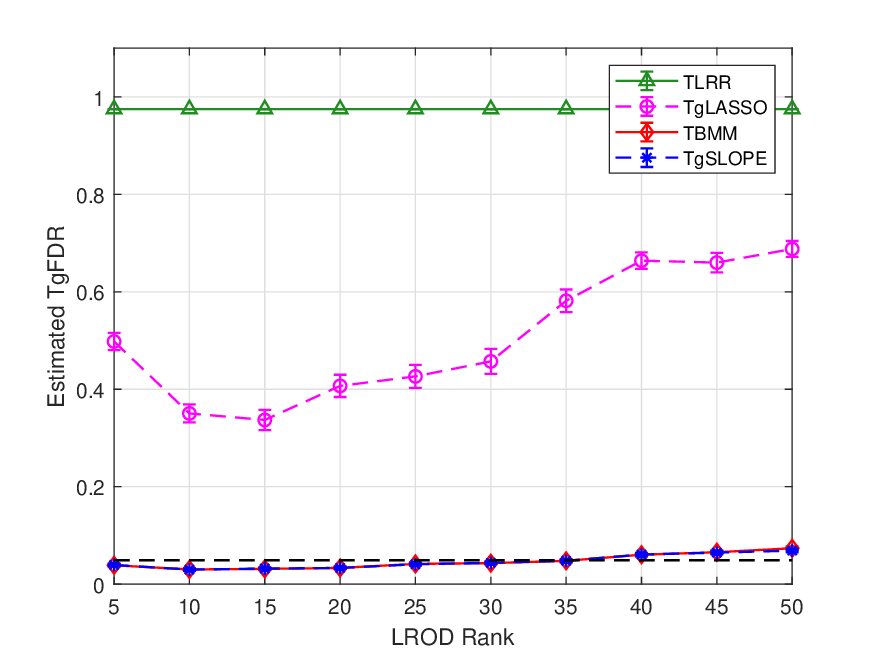}}
          \subfigure[RgEE]{
          \includegraphics[width=2.2in]{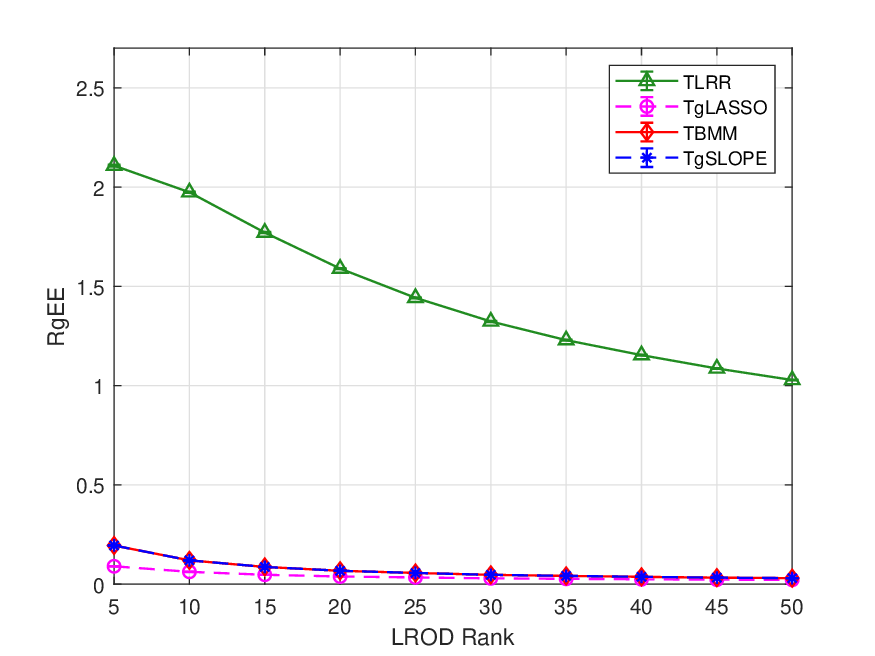}}
          \subfigure[RgEE]{
          \includegraphics[width=2.2in]{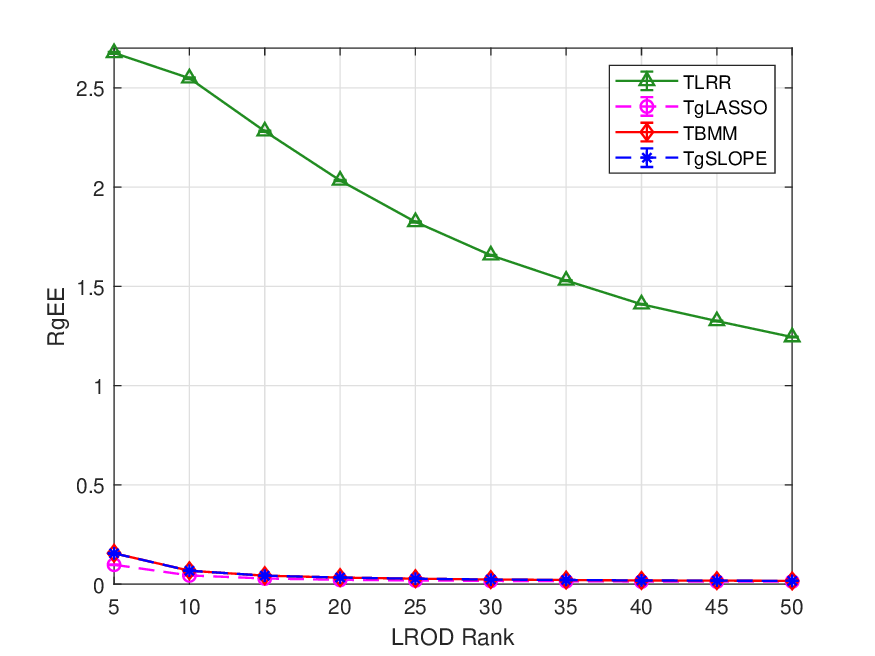}}
          \subfigure[MSE]{
          \includegraphics[width=2.2in]{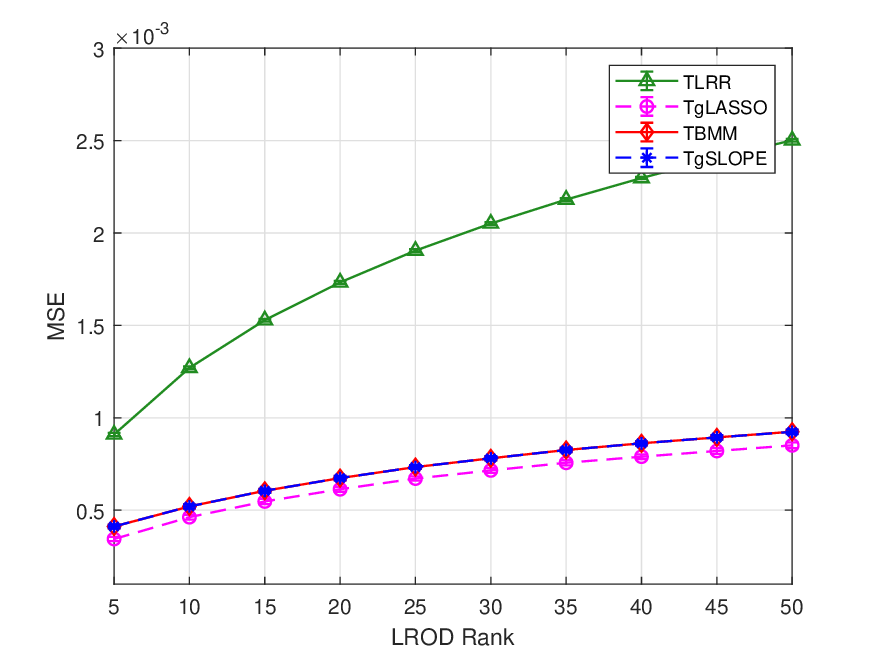}}
          \subfigure[MSE]{
          \includegraphics[width=2.2in]{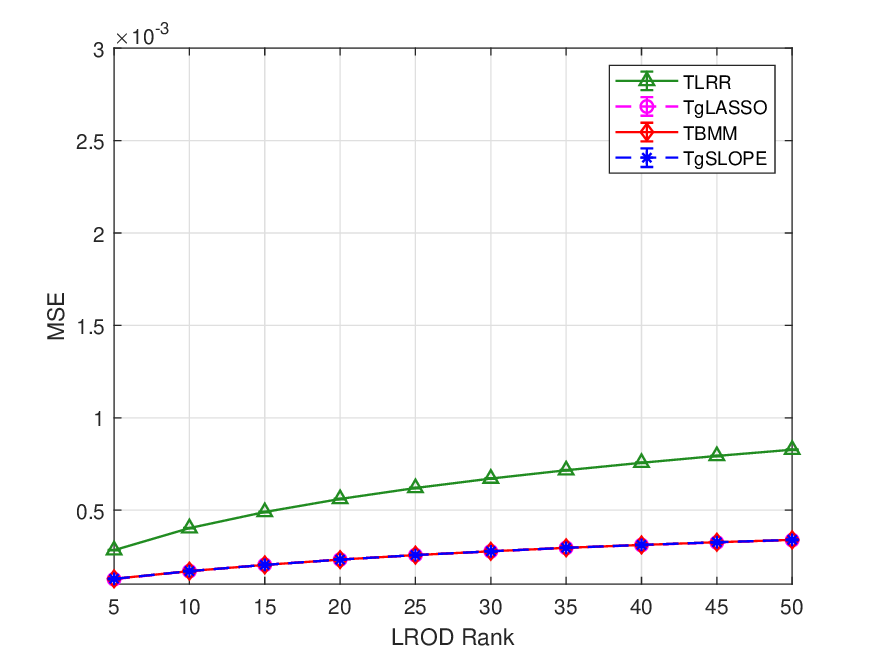}}
          \subfigure[Time]{
          \includegraphics[width=2.2in]{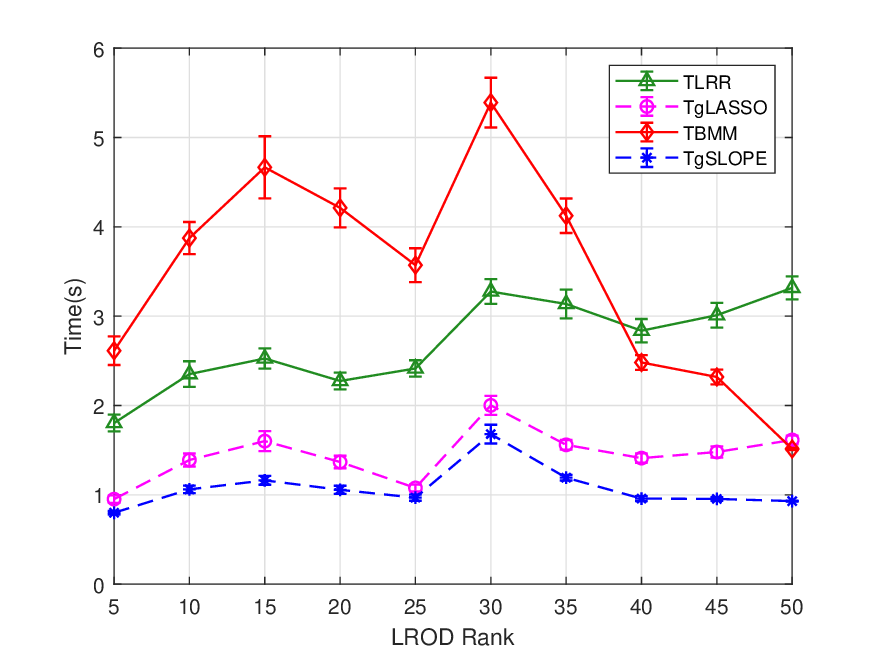}}
          \subfigure[Time]{
          \includegraphics[width=2.2in]{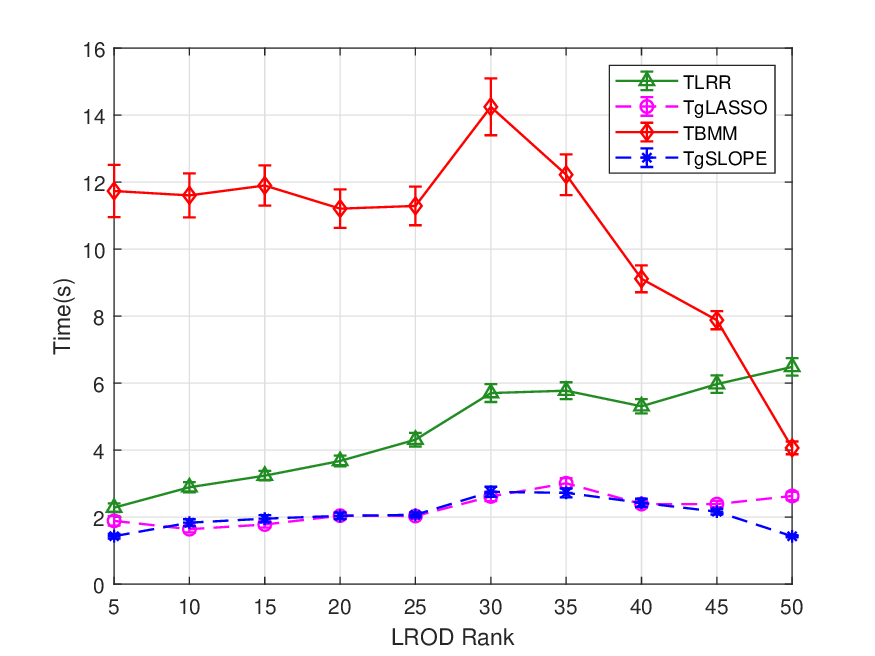}}
           \vspace{-2mm}
           \caption{\small{Simulation results with different LROD rank under $n=1000, p=2000, s=0.02p$ for the Gaussian random design situation of $\bm{X}$. Left column: $p_1=p_2=10$, right column: $p_1=p_2=20$. In (a) and (b), black dotted lines represent the `nominal' TgFDR level $q\cdot(p-s)/p$.}}\label{fig5}
\end{center}
\end{figure}

LROD rank effect: we now examine the impact of LROD rank on the four approaches under $p_1=p_2=10$ and $p_1=p_2=20$, respectively. Set $n=1000, p=2000$ and the sparsity $s=0.02p$. It is known from \cite{2009DecomReview} that LROD rank $K\leq \{p_1p_2, pp_1,pp_2\}$, then we set LROD rank $K=5:5:50$. All testing instances are simulated based on 100 independent replications. Simulation results report that TP =1 for all the competitors. Figure \ref{fig5} depicts changes of the four evaluation metrics with the increase of LROD rank. (i) Our proposed TgSLOPE procedure works well in terms of TgFDR, especially for small LROD ranks. See, e.g., as shown in Figure \ref{fig5} (a) and (b), TgSLOPE and TBMM report the lower TgFDR than the `nominal' level as $K\leq 25$. (ii) Figure \ref{fig5} (g) and (h) illustrate that TgSLOPE reports a relatively short CPU time, especially compared with TBMM. For example, as $p_1=p_2=20$ and $K\leq 35$, the CPU time reported by TgSLOPE is no more than 1/5 of that of TBMM. (iii) For the fixed tensor size, we can see from Figure \ref{fig5} (a) and (b) that the TgFDR tends to become larger as the LROD rank increases. In addition, with the increase of the LROD rank, Figure \ref{fig5} (c) and (d) show that RgEE reduces, while Figure \ref{fig5} (e) and (f) show that MSE grows. This may give us some inspiration to choose a moderate LROD rank.
\subsection{Human Brain Connection Data}
In this subsection, we test our proposed TgSLOPE comparing with the other three approaches on a real human brain connection (HBC) data from the Human Connectome Project (HCP), which aims to build a network map between the anatomical and functional connectivity within healthy human brains \citep{HCP2013}. The preprocessed HBC dataset is provided by \cite{2022JCGS}, in which the response is a $68\times 68$ binary matrix with entries encoding the presence or absence of fiber connections between 68 brain regions-of-interest, the predictor matrix is collected from different personal features for each observed individual. After removing those missing values, the HBC dataset consists of 111 individuals and 549 personal features, including gender, age, etc. (The HBC data can be found at https://wiki.humanconnectome.org/display/PublicData/). \blue{For HBC analysis, the predictor matrix is normalized to have unit column vectors. We choose the tuning parameters via a BIC-type criterion on the whole dataset, which minimizes
$$
{\rm BIC} = \|\Y-\hat{\B}\times_3 \bm{X}\|_F^2+\big(\big\|\|\hat{\B}\|_f\big\|_{0}+p_1+p_2\big)K\log(np_1p_2).
$$
}
In addition, the number of discovered features (\blue{termed as} Discovery), the mean squared prediction error (MSPE) and the CPU time are adopted to evaluate the performance of four competitors. Here, Discovery and MSPE are defined respectively as
\begin{align*}
&{\rm Discovery} = \sharp\{j\in[p]: (\hat{\bm{B}}_{training})_j\neq\bm{O}\}, \\
&{\rm {MSPE}} = \Vert \Y_{test}-\hat{\B}_{training}\times_3\bm{X}_{test}\Vert_F^2/p_1p_2n_{test},
\end{align*}
 where $\hat{\B}_{training}$ is the estimator for the training set, $\bm{X}_{test}$, $\Y_{test}$ and $n_{test}$ are the predictor matrix, response tensor and the sample size of the testing set, respectively.
\begin{figure}
\begin{center}
          \subfigure[MSPE]{
          \includegraphics[width=2.5in]{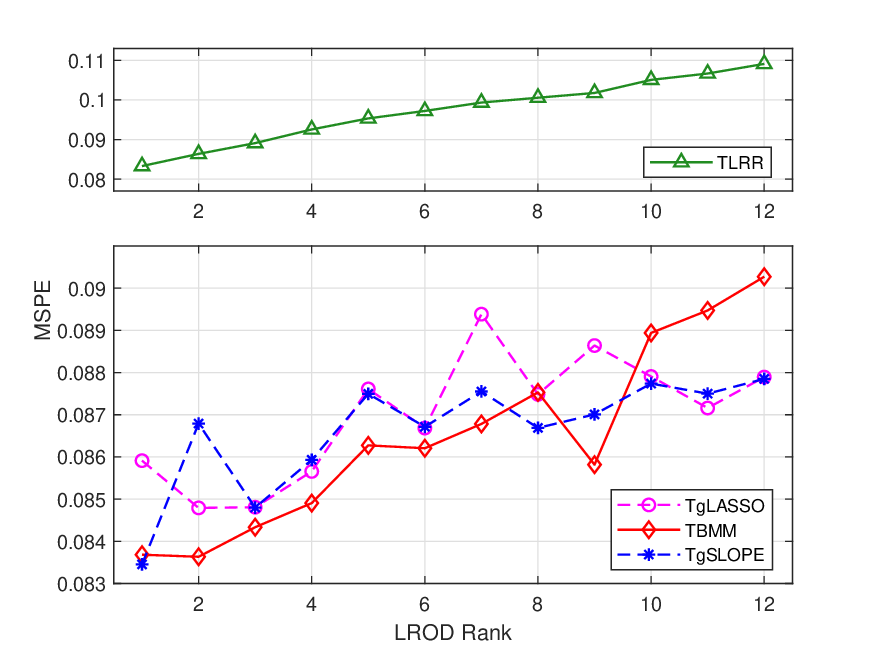}}\hspace{-3mm}
           \subfigure[Discovery]{
          \includegraphics[width=2.5in]{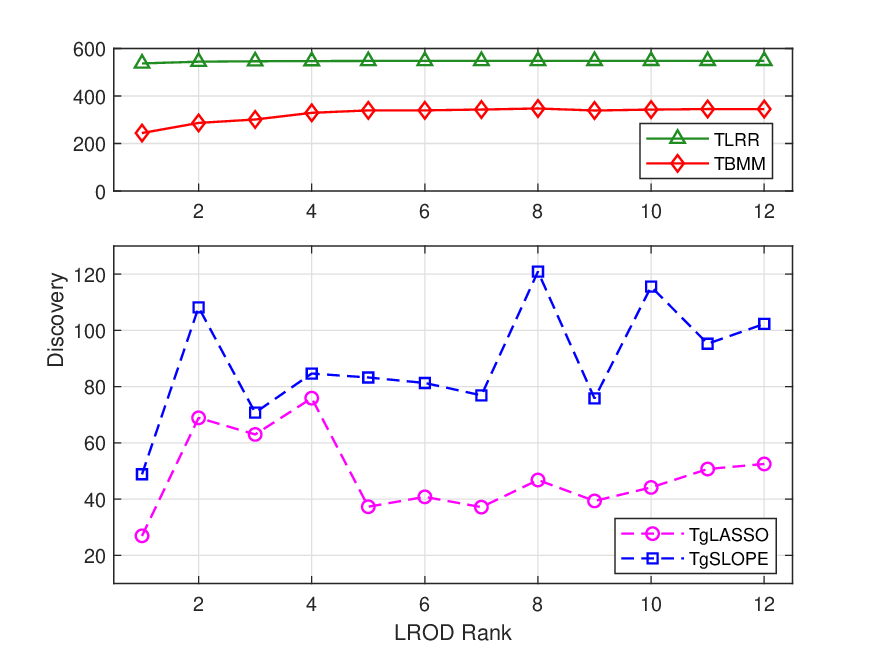}}\hspace{-3mm}
           \caption{\blue{\small{Average results with the different LROD rank based on 20 randomly 8:2 splits for HBC dataset.}}}\label{fig7}
\end{center}
\end{figure}

\blue{For HBC data, the experimental results in \cite{HBC-Sparsity} have revealed that the whole-brain functional magnetic resonance imaging (fMRI) signals can be well characterized by sparse representations, which is supportive on the interpretability with fewer significant features in the HBC data analysis. In this regard, it may be reasonable to recognize the performance superiority of a method with the smaller Discovery and the comparable or even less MSPE. We will test the performance of our proposed approach in such a sense with comparison to other three methods.}

\blue{Notably, the rank of true coefficient tensor is unknown in real data applications.} We first test the performance of four approaches under different LROD ranks. \blue{Set the tolerance for termination in algorithms to be $\epsilon = 10^{-4}$.} For each LROD rank selected from 1 to 12, HBC dataset is randomly divided into 80 percent of training set and 20 percent of testing set 20 times. The average numerical results are depicted in Figure \ref{fig7}. \blue{We can see from Figure \ref{fig7} (a) that small LROD ranks would be sufficient for an acceptable prediction errors of all methods, which suggests us to fit the true rank with a relatively small value. In addition, as shown in Figure \ref{fig7} (b), TLRR fails to feature selection for all rank testing instances, while the other three methods all report the sparse numerical approximate estimates reflected by relatively small values of Discovery, with the smaller MSPE compared to TLRR. This further supports the reasonability of taking the group sparsity of coefficient tensor into consideration for HBC dataset.}

\begin{table}
\renewcommand{\arraystretch}{1.2}
\caption{\blue{Average results with BIC selected LROD rank and different tolerance $\epsilon$ based on 20 randomly 8:2 splits for HBC dataset.}\label{tab-2}}
\setlength{\tabcolsep}{4mm}
\begin{center}
{\footnotesize \begin{tabular}{lllllll}
\hline
    $\epsilon$ &Method      &MSPE        &Discovery    &Time(s)   \\ \hline

    $10^{-4}$  &TLRR        &0.0837      &536.30       &2.766      \\
               &TgLASSO     &0.0848      &28.45        &9.396      \\
               &TBMM        &0.0840      &242.50       &7.137      \\
               &TgSLOPE     &0.0838      &47.95        &8.171      \\  \hline

     $10^{-5}$ &TLRR        &0.0831       &536.20       &9.363     \\
               &TgLASSO     &0.0844       &27.70        &30.561    \\
               &TBMM        &0.0833       &121.90       &92.278    \\
               &TgSLOPE     &0.0833       &47.15        &29.724    \\  \hline

     $10^{-6}$ &TLRR        &0.0836       &537.20       &25.629    \\
               &TgLASSO     &0.0861       &26.70        &111.271   \\
               &TBMM        &0.0838       &73.65        &668.040   \\
               &TgSLOPE     &0.0838       &49.60        &99.964    \\
\hline
 \end{tabular}}
\end{center}
\end{table}

\blue{To test the comparative performance of four methods,} we next choose the combination of LROD rank and regularization parameters via the BIC criterion. \blue{Table \ref{tab-2} collects the average results with different tolerances $\epsilon = \{10^{-4}, 10^{-5},10^{-6}\}$ based on 20 independent replications, which demonstrate the superior performance of TgSLOPE among all comparative methods. Specifically, as shown in Table \ref{tab-2}, TLRR reports the competitive MSPE with TgSLOPE, but unsurprisingly fails to feature selection for all tolerance settings. In addition, TgLASSO gives the larger MSPE than TgSLOPE in all testing instances, although it has the fewest discovered features.

For the method of TBMM, it solves the same optimization model (\ref{cp2-gslope}) by the alternating minimization scheme. While our proposed TgSLOPE method solves problem (\ref{cp2-gslope}) by pDCAe based on the DC reformulation (\ref{dcpg-1}). We can see from Table \ref{tab-2} that these two algorithms report the similar MSPE in different tolerance settings. While the Discovery generated by TBMM is larger than that of TgSLOPE, and the difference of the values for Discovery decreases as the tolerance $\epsilon$ of the numerical algorithms reduces, as shown in Table 2 Columns  ``$\epsilon$" and ``Discovery". In addition, TgSLOPE reports the less CPU time than TBMM, especially for the small tolerance $\epsilon=10^{-6}$. One possible explanation for the superior efficiency of TgSLOPE is that pDCAe has more power of achieving numerical approximate coefficient tensor solution with higher estimate accuracy, leading to better performance on time efficiency and group sparsity reflected by Discovery. The comparison results on Discovery, MSPE and CPU Time of both algorithms with varying $\epsilon$ have been collected in the boxplots in Figure 8.
}
\begin{figure}
\begin{center}
          \includegraphics[width=6in]{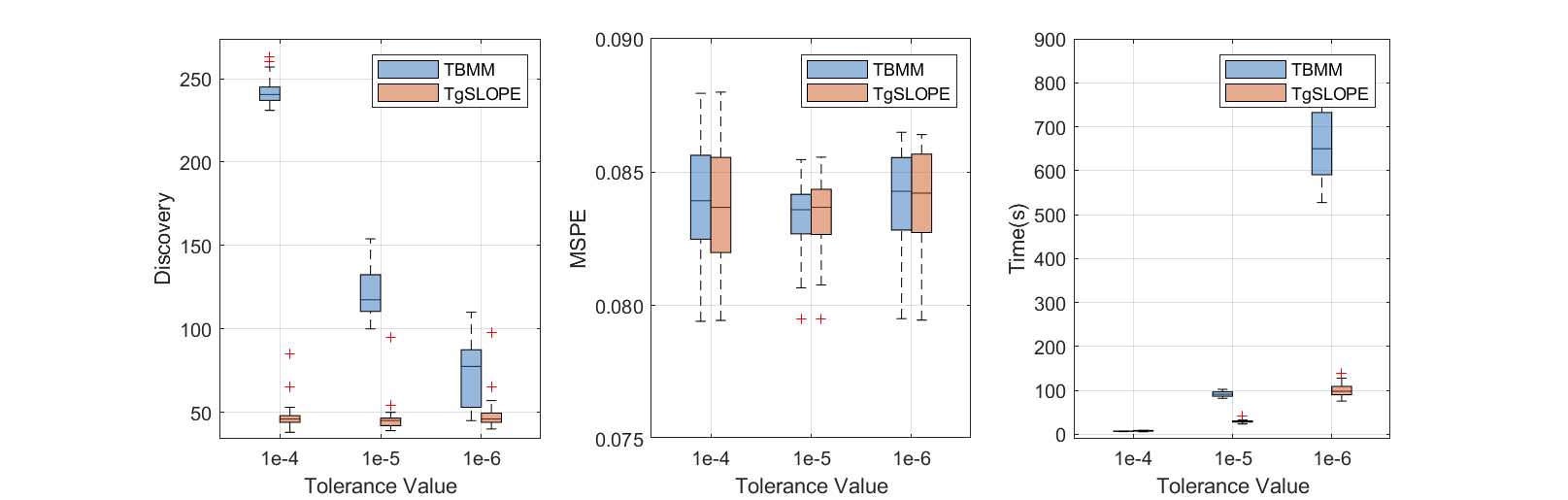}
           \caption{\blue{\small{Boxplot for results of TBMM and TgSLOPE with BIC selected LROD rank and different tolerance $\epsilon$ based on 20 randomly 8:2 splits for HBC dataset.}}}\label{fig8}
\end{center}
\end{figure}
\section{Conclusions \label{sec:6}}
In this article we propose a sparse and low-rank tensor regression method, which optimizes a gSLOPE penalized low-rank, orthogonally decomposable tensor minimization problem. Under the assumption of column-orthogonality for the predictor matrix, we show that our proposed TgSLOPE procedure controls TgFDR at a pre-set level, and achieves the asymptotically minimax convergence with respect to the produced estimation risk. This provides theoretical guarantees for feature selection and coefficient estimation in finite samples. Moreover, a globally convergent pDCAe algorithm is applied to solve the TgSLOPE estimator by constructively reformulating our TgSLOPE problem into a DC program. Numerical experiments verify the superiority of our method in terms of TgFDR control, estimation accuracy and CPU time against three state-of-the-art approaches.

For the TgSLOPE method, it would be interesting to investigate statistical properties including TgFDR control and estimate accuracy in more general cases besides the orthogonal design. Moreover, \cite{2019slopeJMLR} propose a sparse semismooth Newton-based augmented Lagrangian method (Newt-ALM) to solve the SLOPE model of the classical linear regression \citep{2015slope} and show that Newt-ALM offers a notable computational advantage in the high-dimensional settings comparing with the first-order algorithms, such as APG and alternating direction method of multipliers (ADMM). How to design the robust and highly efficient Newton-type algorithm for TgSLOPE based models is also one of our research topic in the future.


\acks{The authors are very grateful to the editor and two anonymous reviewers for their constructive comments and suggestions which greatly help us to improve the quality of our paper. This research was supported by 
the National Natural Science Foundation of China (Grant No.12271022). 
}



\appendix \blue{
\section*{Appendix A. Proof of Proposition \ref{prop-identif}} \label{A-A}
\label{appa:prop}
To check the identifiability for the frontal slice sparse and LROD tensor coefficient $\B\in \mathbb{B}_T$ defined in (\ref{sp-lrtr}), we need to state that the equation system produced by tensor regression model (\ref{tm-2}), i.e., $\bar{\Y} = \B\times_3\bm{X}$,
has a unique solution $\B^*$ satisfying $\B^*\in \mathbb{B}_T$. Under the uniqueness guarantee k$(\bm{W}^*)\geq 2$ of the LROD tensor decomposition $\B^*= [\![\bm{U}^*,\bm{V}^*,\bm{W}^*]\!]$ \cite[Theorem 4b]{UCP-1977}, it suffices to show that
$$\B^* =\argmin_{\B\in\mathbb{R}^{p_1\times p_2\times p}} \big\|\|\B\|_f\big\|_0, \quad {\rm s.t.}  \  \bar{\Y}=\B\times_3\bm{X}$$
is a unique solution. Following the matricized form of a tensor, the above optimization problem can be equivalent to
\begin{align}\label{A-sparsity-recovery}
\min_{\bm{B}\in \mathbb{B}_M} \big\|\|\bm{B}\|_r\big\|_0,
\end{align}
where $\mathbb{B}_M = \{\B\in\mathbb{R}^{p\times p_1p_2}:\M_3(\bar{\Y})=\bm{X}\bm{B}\}$. Thus, it suffices to prove that problem (\ref{A-sparsity-recovery}) has a unique solution $\bm{B}^*=\M_3^{-1}(\B^*)$ under the condition that $\big\|\|\bm{B}^*\|_r\big\|_0\leq k/2$, where $k$ is the k-rank of the matrix $\bm{X}$. 

Assume on the contrary that $\tilde{\bm{B}}\in \mathbb{B}_M$ is another optimal solution of (\ref{A-sparsity-recovery}), which gives that $\big\|\|\bm{B}^*\|_r\big\|_0=\big\|\|\tilde{\bm{B}}\|_r\big\|_0$, $\bm{\Delta} = \bm{B}^*-\tilde{\bm{B}} \neq \bm{O}$ and $\bm{X}\bm{\Delta}=\bm{O}$. This further indicates that there exist $\big\|\|\bm{\Delta}\|_r\big\|_0$ columns of $\bm{X}$ which are linearly dependent. By virtue of the definition of the k-rank of $\bm{X}$, it yields that
\begin{align}\label{A-sparsity-condition}
k  <\big\|\|\bm{\Delta}\|_r\big\|_0 
  \leq \big\|\|\bm{B}^*\|_r\big\|_0+\big\|\|\tilde{\bm{B}}\|_r\big\|_0 \leq k,
\end{align}
where the second inequality is from the fact that $\ell_0$-norm obeys the triangle inequality, and the last inequality from  $\big\|\|\tilde{\bm{B}}\|_r\big\|_0=\big\|\|\bm{B}^*\|_r\big\|_0\leq k/2$.  Contradiction arrives and the proof is completed. 

\hfill$\Box$}

\section*{Appendix B. Proofs for Section \ref{sec:3}} \label{A-B}
\label{appb:theorem}
\subsection*{B.1 Proof of Theorem \ref{th-tfdr}}
Assume that $(\hat{\bm{W}}, \hat{\bm{H}})$ is a local minimizer of the TgSLOPE problem (\ref{cp2-gslope}). Let $\hat{\bm{H}}^{\perp}\in\mathbb{R}^{p_1p_2 \times (p_1p_2-K)}$ such that $[\hat{\bm{H}}, \hat{\bm{H}}^{\perp}]$ is orthogonal. Then the loss function of problem (\ref{cp2-gslope})
\begin{align*}
\begin{split}
L(\bm{W},\hat{\bm{H}})&=\big\Vert[\M_3(\Y)-\bm{X}\bm{W}\hat{\bm{H}}^\top][\hat{\bm{H}}, \hat{\bm{H}}^{\perp}]\big\Vert_F^2\\
& = \big\Vert\M_3(\Y)\hat{\bm{H}}-\bm{X}\bm{W}\big\Vert_F^2+\big\Vert\M_3(\Y)\hat{\bm{H}}^{\perp}\big\Vert_F^2.
\end{split}
\end{align*}
Together with the column-orthogonality of the predictor matrix $\bm{X}$, we have
\begin{align}\label{optg-1}
\hat{\bm{W}}=\argmin_{\bm{W}\in\mathbb{R}^{p\times K}}\bigg\{\frac{1}{2}\big\Vert\hat{\bm{Y}}-\bm{W}\big\Vert_F^2
+P_{\bm{\lambda}}\big(\|\bm{W}\|_r\big)\bigg\},
\end{align}
where $\hat{\bm{Y}}=\bm{X}^\top\M_3(\Y)\hat{\bm{H}}\in\mathbb{R}^{p\times K}$. It follows from (\ref{tm-2}) that $$\hat{\bm{Y}} = \bm{W}^*\bm{H}^{*\top}\hat{\bm{H}}+\bm{X}^\top\M_3(\E)\hat{\bm{H}},$$ which implies that $\hat{\bm{Y}}$ follows the matrix normal distribution $\N(\bm{W}^*\bm{H}^{*\top}\hat{\bm{H}}, \sigma^2\bm{I}_p\otimes\bm{I}_K)$. Note that (\ref{optg-1}) is a convex problem with strongly convex, piecewise linear-quadratic objective function and hence admits a unique optimal solution. Similar to the proximal operator for gSLOPE \citep{2019gslope}, the optimization problem (\ref{optg-1}) can be solved in two steps
\begin{align}\label{optg-2}
\left\{
\begin{array}{l}
\hat{\bm{\eta}}=\argmin_{\bm{\eta}}\big\{\frac{1}{2}\sum_{j=1}^p \big( \Vert\hat{\bm{y}}_j\Vert-\eta_j\big)^2+P_{\bm{\lambda}}(\bm{\eta})\big\},\vspace{1ex}\\
\hat{\bm{W}}_j=\hat{\eta}_j\frac{\hat{\bm{y}}_j}{\Vert\hat{\bm{y}}_j\Vert}, j\in[p].
\end{array}
\right.
\end{align}

Define $V_{\bm{\eta}}=\sharp\{j\in[p]: \bm{W}_j^*=\bm{0}, \hat{\eta}_j\neq0\}, R_{\bm{\eta}}=\sharp\{j\in[p]: \hat{\eta}_j\neq0\}$. It suffices to show that
$$\mathbb{E}\bigg[\frac{V_{\bm{\eta}}}{\max\{R_{\bm{\eta}},1\}}\bigg]\leq q\cdot \frac{p-s}{p}.$$
Without loss of generality, assume that $\bm{W}_j^*=\bm{0}$ for $0\leq j\leq p-s$ and $\bm{W}_j^*\neq \bm{0}$ otherwise. By the definition of TgFDR, we derive that
\begin{align}\label{tFDR}
\begin{split}
\mathbb{E}\bigg[\frac{V_{\bm{\eta}}}{\max\{R_{\bm{\eta}},1\}}\bigg] &= \sum_{r=1}^{p}\frac{1}{r}\mathbb{E}\big[V_{\bm{\eta}}\mathbbm{1}_{\{R_{\bm{\eta}} =r\}}\big]\\
& = \sum_{r=1}^{p}\frac{1}{r}\sum_{j=1}^{p-s}\mathbb{E}\big[\mathbbm{1}_{\{\hat{\eta}_j\neq 0\}}\mathbbm{1}_{\{R_{\bm{\eta}} =r\}}\big]\\
& = \sum_{r=1}^{p}\frac{1}{r}\sum_{j=1}^{p-s}\Pr\big(\hat{\eta}_j\neq 0, R_{\bm{\eta}} =r\big),
\end{split}
\end{align}
where $\mathbbm{1}_{\{\cdot\}} =1$ if the event occurs, and $\mathbbm{1}_{\{\cdot\}}=0$ otherwise. Next we focus on the events $\{\hat{\eta}_j\neq 0, R_{\bm{\eta}} =r\}$, $j=1,\dots,p-s, r=1,\dots,p$. Denote $\overline{\bm{Y}}=(\hat{\bm{y}}_1, \dots,\hat{\bm{y}}_{j-1},\hat{\bm{y}}_{j+1},\dots,\hat{\bm{y}}_p)^\top\in\mathbb{R}^{(p-1)\times K}$, $\overline{\bm{\lambda}}=(\lambda_2, \dots, \lambda_p)^\top\in\mathbb{R}^{p-1}$. Applying the simplified TgSLOPE problem (\ref{optg-1}) to $\overline{\bm{Y}}$ with tuning parameter $\overline{\bm{\lambda}}$, the optimization problem is given by
\begin{align*}
\overline{\bm{W}}=\argmin_{\bm{W}\in\mathbb{R}^{(p-1)\times K}} \Bigg\{\frac{1}{2}\sum_{j=1}^{p-1}\big\Vert\overline{\bm{y}}_{j}-\bm{W}_j\big\Vert^2
+P_{\overline{\bm{\lambda}}}\big(\|\bm{W}\|_r\big) \Bigg\}.
\end{align*}
Define $\overline{R}^j = \sharp\{j\in[p-1]: \overline{\bm{W}}_j\neq \bm{0}\}$. It follows from Lemmas E.6 and E.7 in \cite{2019gslope} that
 \begin{align*}
 \big\{\|\hat{\bm{Y}}\|_r: \hat{\eta}_j\neq 0, R_{\bm{\eta}} = r\big\} \subset \big\{\|\hat{\bm{Y}}\|_r: \Vert\hat{\bm{y}}_j\|>\lambda_{r}, \overline{R}^j =r-1\big\},
 \end{align*}
 where $\|\hat{\bm{Y}}\|_r = (\Vert\hat{\bm{y}}_1\Vert, \dots,\Vert\hat{\bm{y}}_p\Vert)^\top$ with $\Vert\hat{\bm{y}}_j\Vert^2\sim \chi^2_K\big(\Vert\bm{W}_j^*\bm{H}^{*\top}\hat{\bm{H}}\Vert^2\big)$, $j=1,\dots,p$. Then we have
 \begin{align*}
\Pr(\hat{\eta}_j\neq 0, R_{\bm{\eta}} = r)&\leq \Pr(\Vert\hat{\bm{y}}_j\Vert>\lambda_{r}, \overline{R}^j =r-1)\\
&= \Pr(\Vert\hat{\bm{y}}_j\Vert/\sigma>\lambda_{r}^*)\Pr(\overline{R}^j =r-1)\\
&\leq \frac{q\cdot r}{p}\Pr(\overline{R}^j =r-1),
 \end{align*}
 where $\lambda_{r}^* = F_{\chi_{K}}^{-1}\big(1-q\cdot r/p\big)$, the equality is due to the independence between $\Vert\hat{\bm{y}}_j\Vert$ and $\overline{R}^j$, the last inequality follows from the definition of the tuning parameters in (\ref{tuning}). Therefore, $\TgFDR$ in (\ref{tFDR}) can be bounded by
  \begin{align*}
\mathbb{E}\bigg[\frac{V_{\bm{\eta}}}{\max\{R_{\bm{\eta}},1\}}\bigg] &\leq \sum_{r=1}^{p}\frac{1}{r}\sum_{j=1}^{p-s}\frac{q\cdot r}{p}\Pr(\overline{R}^j =r-1)\\
 &= \sum_{j=1}^{p-s}\frac{q}{p}\sum_{r=1}^{p}\Pr(\overline{R}^j =r-1)=q\cdot \frac{p-s}{p}.
 \end{align*}
The proof is completed.  \hfill$\Box$

\subsection*{B.2 Proof of Theorem \ref{thm-2}}
Let $(\hat{\bm{W}}, \hat{\bm{H}})$ be an optimal solution of the problem (\ref{cp2-gslope}). Then $\hat{\B}= \M_3^{-1}(\hat{\bm{W}}\hat{\bm{H}}^\top)$. Under the column-orthogonality assumption on the predictor matrix $\bm{X}$, we know that the statistically equivalent model of (\ref{tm-2}) is
\begin{align}\label{eq1}
\widetilde{\bm{Y}} =\M_3(\B^*)+\bm{X}^\top\M_3(\E),
\end{align}
which has the distribution $\N(\M(\B^*), \sigma^2\bm{I}_p\otimes\bm{I}_K)$. Denote $\widetilde{\mathbb{B}}$ as the set of all coefficient tensors for which only the elements in the first column of the mode-3 unfolded matrix are possibly nonzero and the rest are fixed to be zero. Let $\widetilde{\mathbb{B}}_{s} = \mathbb{B}_{s}\cap\widetilde{\mathbb{B}}$. Then, for any $\B^*\in \widetilde{\mathbb{B}}_{s}$, (\ref{eq1}) is reduced to a general Gaussian sequence model with length $p$ and sparsity at most $s$. As $s/p\to 0$, this sequence model has minimax risk $(1+o(1))2\sigma^2s\log(p/s)$ \citep{1994minimax}. Thus, we have
$$\sup_{\B^*\in \widetilde{\mathbb{B}}_{s}}\mathbb{E}\big\Vert \M(\hat{\B})-\M(\B^*) \big\Vert_F^2 \blue{\sim} (1+o(1))2\sigma^2s\log(p/s),$$
which yields that
$$\sup_{\B^*\in \mathbb{B}_{s}}\mathbb{E}\blue{\|\hat{\B}-\B^*\|_F^2} \geq (1+o(1))2\sigma^2s\log(p/s).$$

We next show that the $\ell_2$-loss which measures the deviation of the TgSLOPE estimator from the ground truth $\B^*$ is bounded above by $(1+o(1))2\sigma^2s\log(p/s)$. For simplicity, we assume that $\|\bm{w}_j^*\|\neq 0$ for $j\leq s$ and $\|\bm{w}_j^*\|=0$ otherwise. Denote \blue{$\mu_j=\|\bm{w}_j^*\|$} and \blue{$\nu_j=\|\hat{\bm{H}}^\top\bm{H}^*\bm{w}_j^*\|$}. Then, it follows from the proof of Theorem \ref{th-tfdr} that $\hat{\bm{W}}$ can be obtained by the two step format (\ref{optg-2}), in which $\Vert\hat{\bm{y}}_j\Vert^2\sim \chi^2_K(\nu_j^2)$, $j=1,\dots,p$. The $\ell_2$-loss is
\blue{\begin{align*}
\mathbb{E}\|\hat{\B}-\B^*\|_F^2&\overset{(a)}{=} \mathbb{E}\|\hat{\bm{W}}\hat{\bm{H}}^\top-\bm{W}^*\bm{H}^{*\top}\|_F^2
=\mathbb{E}\Bigg[\sum_{j=1}^p\|\hat{\bm{H}}\hat{\bm{w}}_j-\bm{H}^*\bm{w}_j^*\|^2\Bigg]\\
&\overset{(b)}{\leq}\mathbb{E}\Bigg[\sum_{j=1}^p(\|\hat{\bm{H}}\hat{\bm{w}}_j\|+\|\bm{H}^*\bm{w}_j^*\|)^2\Bigg]
\overset{(c)}{=}\mathbb{E}\Bigg[\sum_{j=1}^p(\|\hat{\bm{w}}_j\|+\|\bm{w}_j^*\|)^2\Bigg],\\
&=\mathbb{E}\Bigg[\sum_{j=1}^{s}(\|\hat{\bm{w}}_j\|+\|\bm{w}_j^*\|)^2\Bigg]+\mathbb{E}\Bigg[\sum_{j=s+1}^{p}\|\hat{\bm{w}}_j\|^2\Bigg]\\
&=\underbrace{\mathbb{E}\Bigg[\sum_{j=1}^{s}(\hat{\eta}_j+\mu_j)^2\Bigg]}_{E1}+\underbrace{\mathbb{E}\Bigg[\sum_{j=s+1}^{p}\hat{\eta}_j^2\Bigg]}_{E2},
\end{align*}
where (a) comes from $\hat{\B}= \M_3^{-1}(\hat{\bm{W}}\hat{\bm{H}}^\top)$, (b) is due to the triangle inequality and (c) follows from the column-orthogonality of $\hat{\bm{H}}$ and $\bm{H}$.} Thus, it suffices to show
$$E1\leq (1+o(1))2\sigma^2s\log(p/s) \ \ and \ \ E2= o(1)2\sigma^2s\log(p/s).$$

To proceed, define random variables $\psi_j^2=\psi_{j,1}^2+\psi_{j,2}^2+\cdots+\psi_{j,K}^2$ and $\phi_j^2=(\psi_{j,1}+\nu_j)^2+\psi_{j,2}^2+\cdots+\psi_{j,K}^2$ with i.i.d. $\psi_{j,k}\sim N(0,1), j\in[p], k\in[K]$. Then $\psi_j^2\sim \chi_K^2$ and $\phi_j^2\sim \chi_K^2(\nu_j^2)$ for all $j\in[p]$. Denoting $\bm{\phi} = (\phi_1, \phi_2, \dots, \phi_p)^\top$, we have
\begin{align}\label{E1-1}
\begin{split} \mathbb{E}\big\|\hat{\bm{\eta}}_{[1:s]}\blue{+}\bm{\mu}_{[1:s]}\big\|^2&=\mathbb{E}\big\|\hat{\bm{\eta}}_{[1:s]}-\bm{\phi}_{[1:s]}+\bm{\phi}_{[1:s]}\blue{+}\bm{\mu}_{[1:s]}\big\|^2\\
&\leq \mathbb{E}\Big(\big\|\hat{\bm{\eta}}_{[1:s]}-\bm{\phi}_{[1:s]}\big\|+\big\|\bm{\phi}_{[1:s]}\blue{+}\bm{\mu}_{[1:s]}\big\|\Big)^2\\
&\leq \big\|\bm{\lambda}_{[1:s]}\big\|^2+\mathbb{E}\big\|\bm{\phi}_{[1:s]}\blue{+}\bm{\mu}_{[1:s]}\big\|^2+2\big\|\bm{\lambda}_{[1:s]}\big\|\mathbb{E}\big\|\bm{\phi}_{[1:s]}\blue{+}\bm{\mu}_{[1:s]}\big\|,
\end{split}
\end{align}
where the last inequality is obtained by $\big\|\hat{\bm{\eta}}_{[1:s]}-\bm{\phi}_{[1:s]}\big\|\leq \big\|\bm{\lambda}_{[1:s]}\big\|$, owing to Fact 3.3 of \cite{2016slope} with conditions $\phi_j^2$ and $\Vert\hat{\bm{y}}_j\Vert^2$ are i.i.d. for any $j\leq s$.  Moreover, it follows from \cite{2010chi} that as $s/p\to 0$, $F_{\chi_{K}}^{-1}\big(1-q\cdot j/p\big)\sim \sqrt{2\log(p/(q\cdot j))}$ for all $j\leq s$, which yields
\begin{align}\label{proof-0}
\big\|\bm{\lambda}_{[1:s]}\big\|^2 \sim 2\sigma^2s\log(p/s).
\end{align}
Then, it is easy to see that
\begin{align}\label{proof-1}
\begin{split}
|\phi_j\blue{+}\mu_j|^2&=\Big|\sqrt{(\psi_{j,1}+\nu_j)^2+\psi_{j,2}^2+\cdots+\psi_{j,K}^2}\blue{+}\mu_j\Big|^2\\
&\leq \Big|\sqrt{\psi_{j,2}^2+\cdots+\psi_{j,K}^2}+|\psi_{j,1}|+\nu_j\blue{+}\mu_j\Big|^2\\
&\leq \Big(\sqrt{2(\psi_{j,1}^2+\psi_{j,2}^2+\cdots+\psi_{j,K}^2)}+\blue{\nu_j+\mu_j}\Big)^2\\
& \leq 4\psi_{j}^2+2(\nu_j\blue{+}\mu_j)^2.
\end{split}
\end{align}
Plugging \blue{$\mu_j=\|\bm{w}_j^*\|$} and \blue{$\nu_j=\|\hat{\bm{H}}^\top\bm{H}^*\bm{w}_j^*\|$} into the part $|\nu_j-\mu_j|$ derives that
\begin{align}\label{proof-2}
\begin{split}
|\nu_j\blue{+}\mu_j|
&\overset{(a)}{\leq} \blue{\Vert\hat{\bm{H}}^\top\bm{H}^*\bm{w}_j^*+\bm{w}_j^*\Vert
=\big\|(\hat{\bm{H}}+\bm{H}^*)^\top(\M_3(\B^*))_j\big\|}\\
&\leq \big\|(\M_3(\B^*))_j\big\|\Vert\hat{\bm{H}}\blue{+}\bm{H}^*\Vert_2\leq \big\|(\M_3(\B^*))_j\big\|\big(\Vert \hat{\bm{H}}\Vert_2+\Vert \bm{H}^*\Vert_2\big)\\
& \overset{(b)}{\leq} 2\big\|(\M_3(\B^*))_j\big\|,
\end{split}
\end{align}
where (a) follows from the Cauchy-Schwarz inequality, (b) comes from the fact that $\Vert \hat{\bm{H}}\Vert_2=\Vert \bm{H}^*\Vert_2=1$ for the column-orthogonal matrices $\hat{\bm{H}}$ and $\bm{H}^*$. Setting $\bar{\alpha} = \max\{\Vert\bm{B}_j^*\Vert_F, j\in[p]\}$, we know that $(\nu_j-\mu_j)^2\leq 4\bar{\alpha}^2$. Combining with (\ref{proof-1}) and (\ref{proof-2}), we obtain
\begin{align}\label{proof-3}
\begin{split}
\mathbb{E}\big\|\bm{\phi}_{[1:s]}\blue{+}\bm{\mu}_{[1:s]}\big\|^2&= \mathbb{E}\bigg[\sum_{j=1}^s(\phi_j\blue{+}\mu_j)^2\bigg]\leq \mathbb{E}\bigg[\sum_{j=1}^s(4\psi_j^2+8\bar{\alpha}^2)\bigg]\\
&= 8s\bar{\alpha}^2+4\mathbb{E}(\zeta^2) =4s(2\bar{\alpha}^2+K),
\end{split}
\end{align}
where $\zeta^2=\sum_{j=1}^s\psi_j^2\sim\chi^2_{sK}$.
Therefore, combining with (\ref{E1-1}), (\ref{proof-0}) and (\ref{proof-3}) yields that
\begin{align}\label{E1-2}
\begin{split}
E1&\leq (1+o(1))2\sigma^2s\log(p/s)+4s(2\bar{\alpha}^2+K)\\
&+2\sqrt{4s(2\bar{\alpha}^2+K)}\sqrt{(1+o(1))2\sigma^2s\log(p/s)}\\
&\sim (1+o(1))2\sigma^2s\log(p/s),
\end{split}
\end{align}
where the last step makes use of $s/p\to 0$.

We claim that $E2= o(1)2\sigma^2s\log(p/s)$ in the following. Note that $|\phi_j|=|\psi_j|\sim \chi_K$ since $\nu_j=0$ for $j>s$. Denote $|\psi|_{(1)}\geq\cdots\geq |\psi|_{(p-s)}$ as the order statistics of $|\psi_{s+1}|,\dots,|\psi_{p}|$. It follows from the proof of Lemma 3.3 in \cite{2016slope} that
$$\mathbb{E}\Bigg[\sum_{j=s+1}^{p}\hat{\eta}_j^2\Bigg]\leq \sum_{j=1}^{p-s}\mathbb{E}(|\psi|_{(j)}-\lambda_{s+j})_{+}^2,$$
where $x_{+}=\max\{0,x\}$. Then, we can partition the sum into three parts
\begin{align*}
\sum_{j=1}^{p-s}\mathbb{E}(|\psi|_{(j)}-\lambda_{s+j})_{+}^2&=\sum_{j=1}^{\lfloor As\rfloor}\mathbb{E}(|\psi|_{(j)}-\lambda_{s+j})_{+}^2\\
&+\sum_{j=\lceil As\rceil }^{\lfloor ap\rfloor }\mathbb{E}(|\psi|_{(j)}-\lambda_{s+j})_{+}^2+\sum_{j=\lceil ap\rceil }^{p-s}\mathbb{E}(|\psi|_{(j)}-\lambda_{s+j})_{+}^2,
\end{align*}
for a sufficiently large constant $A>0$ and a sufficiently small constant $a>0$. Note that Lemmas F.1, F.2 and F.3 given by \cite{2019gslope} show that each part is negligible compared with $2\sigma^2s\log(p/s)$. This indicates that $E2= o(1)2\sigma^2s\log(p/s)$ and consequently completes the proof together with (\ref{E1-2}). \hfill$\Box$

\section*{Appendix C. Proofs for Section \ref{sec:4}}\label{A-C}
\label{appc:theorem}
\subsection*{C.1 Proof of Lemma \ref{lemma-3}}
Denoting $\sigma_1\geq\sigma_2\geq \cdots\geq\sigma_{K}$ as singular values of $\bm{W}$, we can obtain that
\begin{align}\label{converg-1}
\|\bm{W}\|_{*}^2=\Bigg(\sum_{i=1}^K\sigma_i\Bigg)^2\leq K \sum_{i=1}^K\sigma_i^2=K\|\bm{W}\|_F^2.
\end{align}
In addition, we know from the singular value inequality \citep{Eigvalues} that $\sigma_{i+j-1}(\bm{A}\bm{B})\leq \sigma_{i}(\bm{A})\sigma_{j}(\bm{B})$ for $1\leq i,j\leq K, i+j\leq K+1$, which implies that $\|\bm{A}\bm{B}\|_{*}\leq \|\bm{A}\|_{*}\|\bm{B}\|_{*}$ for any two matrices $\bm{A}\in\mathbb{R}^{p_1p_2\times p}$ and $\bm{B}\in\mathbb{R}^{p\times K}$. Together with (\ref{converg-1}), we derive that as $\|\bm{W}\|_F\to\infty$,
\begin{align*}
\frac{F(\bm{W})}{\|\bm{W}\|_F}&=\frac{\|\bm{X}\bm{W}\|_F^2/2+P_{\bm{\lambda}}\big(\|\bm{W}\|_r\big)-\|\M_3(\Y)^\top\bm{X}\bm{W}\|_{*}}{\|\bm{W}\|_F}\\
&\geq \frac{\|\bm{X}\bm{W}\|_F^2/2+P_{\bm{\lambda}}\big(\|\bm{W}\|_r\big)}{\|\bm{W}\|_F}-\frac{\|\M_3(\Y)^\top\bm{X}\|_{*}\|\bm{W}\|_{*}}{\|\bm{W}\|_F}\\
&\overset{}{\geq }\frac{\|\bm{X}\bm{W}\|_F^2/2+P_{\bm{\lambda}}\big(\|\bm{W}\|_r\big)}{\|\bm{W}\|_F}-\sqrt{K}\|\M_3(\Y)^\top\bm{X}\|_{*}\to \infty,
\end{align*}
where `$\to$' is from the fact that the inequality $\|\bm{X}\bm{W}\|_F^2\geq \sigma_p^2(\bm{X})\|\bm{W}\|_F^2$ holds with $\sigma_p(\bm{X})>0$ when $\bm{X}$ is full column rank. This completes the proof. \hfill$\Box$

\subsection*{C.2 Proof of Theorem \ref{thm-3}}
It follows from Lemma \ref{lemma-3} that the DC function $F$ in (\ref{pdca-2}) is also level-bounded, that is, the level set $\L_{\alpha}=\{\bm{W}\in\mathbb{R}^{p\times K}: F(\bm{W})\leq \alpha\}$ is bounded for any $\alpha\in {\mathbb{R}}$. In addition, we know from \cite{Chen-pdca} that pDCAe enjoys the global convergence if $F$ is level-bounded. Thus, using Lemma \ref{lemma-3}, the desired properties in (a)-(c) are consequences of Theorem 4.1 in \cite{Chen-pdca}. \hfill $\Box$

\blue{
\section*{Appendix D. Scaling of nonzero ground truth in simulations}\label{A-D}
In simulations, we generate the row sparsity factor matrix $\bm{W}^*\in\mathbb{R}^{p\times K}$ in a similar manner as in \cite{2019gslope}, where each nonzero row of $\bm{W}^*$ is scaled such that $\|\bm{w}_j^*\|=a\sqrt{K}$ with $a=\sqrt{4\ln(p)/(1-p^{-2/K})-K}$. This specific form allows the generated signals to be comparable to the maximal noise such that nonzero signals (i.e., significant variables) can be identified with moderate power. We give the calculation of the scale value $a$ in the following.

Considering the case of the column-orthogonal predictor matrix, the TgSLOPE estimator $\hat{\bm{W}}$ can be obtained from (\ref{optg-1}). We see from (\ref{optg-1}) that the identification of
the $R$ significant variables (where the number, $R$, is determined by $\bm{\lambda}$) corresponds to discovering indices of the $R$ largest values among $\|\hat{\bm{y}}_1\|, \ldots, \|\hat{\bm{y}}_p\|$. Note that $\hat{\bm{y}}_1, \dots, \hat{\bm{y}}_p$ are generated respectively by the random vectors $\hat{Y}_j = (\hat{Y}_{j_1},\dots,\hat{Y}_{j_K})^\top,j =1,\dots,p$, where
\begin{align}\label{R-1}
\hat{Y}_j=\bm{v}_j+E_j\sim \N(\bm{v}_j, \sigma^2\bm{I}_{K})
\end{align}
with $\bm{v}_j=\hat{\bm{H}}^\top\bm{H}^*\bm{w}_j^*$ and i.i.d. random noise vectors $E_j = \hat{\bm{H}}^\top(\M_3(\E))^\top\bm{X}_{:j}\sim \N(0, \sigma^2\bm{I}_{K})$ for $j =1,\dots,p$. Thus $\|\hat{Y}_j\|=\sqrt{\sum_{k=1}^K \hat{Y}_{j_k}^2}$ has a $\chi_K$ distribution with the noncentrality parameter $\|\bm{v}_j\|$ and $\|E_j\|$ has a central $\chi_K$ distribution for $j\in[p]$. Then, the nonzero $\|\bm{v}_j\|$ (or the nonzero $\bm{w}_j^*$) could be perceived as a strong signal and thus be identified by TgSlOPE if with high probability the value $\|\hat{\bm{y}}_j\|$ generated from the noncentral $\chi$ distribution $\hat{Y}_j$ is large compared to the background noise produced by the random disturbance $E_j$ with central $\chi_K$ distributions. Otherwise, the signal could be easily covered by random disturbances. Theorem H.1 of \cite{2019gslope} gives that for independent variables $Z_1,\ldots Z_p$ with $Z_j\sim \chi_K^2, j\in[p]$, one has $$\mathbb{E}\Big(\max_{j\in[p]}\{Z_j\}  \Big)\leq \frac{4\ln(p)}{1-m^{-2/K}}.$$
Thus an idea is to use the quantity $\sqrt{4\ln(p)/(1-m^{-2/K})}$ as the upper bound on the expected value of maximum over $p$ independent $\chi_K$-distributed variables. To investigate the discovery performance of TgSLOPE in simulation, we aim at $\mathbb{E}(\|\hat{\bm{y}}_j\|)=\sqrt{4\ln(p)/(1-m^{-2/K})}$, which yields the setting that $\|\bm{v}_j\|=\sqrt{4\ln(p)/(1-m^{-2/K})-K}$ since $\mathbb{E}(\|\hat{\bm{y}}_j\|)\approx \sqrt{\|\bm{v}_j\|+K}$. Note that $\|\bm{v}_j\|=\|\hat{\bm{H}}^\top\bm{H}^*\bm{w}^*_j\|\leq \|\hat{\bm{H}}\|_2\|\bm{H}^*\|_2\|\bm{w}^*_j\|=\|\bm{w}^*_j\|$ due to the column-orthogonality of $\hat{\bm{H}}$ and $\bm{H}^*$. In such a sense, we use the value $\sqrt{4K\ln(p)/(1-p^{-2/K})-K}$ to scale $\|\bm{w}_j^*\|$ in our simulations. This yields signals comparable to the maximal noise such that nonzero signals can be detected with moderate power.}

\vskip 0.2in
\bibliography{Refs}

\end{document}